\documentclass[12pt,a4paper]{article}
\usepackage{amsmath,amsfonts,amssymb,amsthm}
\usepackage{amscd}

\title{On the noncommutative spin geometry of
the standard Podle\'s sphere
and index computations}

\author{Elmar Wagner}

\date{}

\makeatletter
\renewcommand{\[}{\begin{equation}}
\renewcommand{\]}{\end{equation}}
\newcommand{\ip}[2]{\langle #1,#2 \rangle}
\newcommand{\fett}[1]{{\bf #1}}
\newcommand{\msum}{\mbox{$\sum$}}

\newcommand{\srep}{*-re\-pre\-sen\-ta\-tion}

\newtheorem{thm}{Theorem}[section]
\newtheorem{prop}[thm]{Proposition}
\newtheorem{lem}[thm]{Lemma}
\newtheorem{cor}[thm]{Corollary}

\theoremstyle{definition}
\newtheorem{defn}[thm]{Definition}
\newtheorem{rem}[thm]{Remark}
\newtheorem{ex}[thm]{Example}

\newcommand{\hH}{{\mathcal{H}}}
\newcommand{\cO}{{\mathcal{O}}}
\newcommand{\cK}{{\mathcal{K}}}
\newcommand{\U}{{\mathcal{U}}}   

\newcommand{\cI}{{\mathcal{I}}}

\newcommand{\A}{{\mathcal{A}}}
\newcommand{\B}{{\mathcal{B}}}
\newcommand{\cC}{{\mathcal{C}}}

\newcommand{\cS}{{\mathcal{S}}}

\newcommand{\Z}{{\mathbb{Z}}}
\newcommand{\R}{{\mathbb{R}}}
\newcommand{\C}{{\mathbb{C}}}
\newcommand{\N}{{\mathbb{N}}}

\newcommand{\dif}{\partial}

\newcommand{\ov}{\overline}
\newcommand{\im}{\mathrm{i}}
\newcommand{\trans}{\mathrm{t}}
\newcommand{\Lin}{{\mathrm{span}}}
\newcommand{\M}{\mathrm{M}}

\newcommand{\mL}{\mathrm{L}}

\newcommand{\id}{{\mathrm{id}}}
\newcommand{\ad}{{\mathrm{ad}}}      %adjungierte Wirkung

\newcommand{\ind}{{\mathrm{ind}}}
\newcommand{\qind}{q\mbox{-}{\mathrm{ind}}}

\newcommand{\End}{\mathrm{End}}

\newcommand{\adro}{\mbox{$\mathrm{ad}^\circ_{\mathrm{R}}$}}

\newcommand{\cop}{\mathrm{cop}}
\newcommand{\Hom}{\mathrm{Hom}}
\newcommand{\HC}{\mathrm{HC}}
\newcommand{\HH}{\mathrm{HH}}
\newcommand{\BC}{\mathrm{BC}}
\newcommand{\ch}{\mathrm{ch}^\s} \newcommand{\chid}{\mathrm{ch}^\id}
\newcommand{\cht}{\mathrm{ch}^\theta}
\newcommand{\chn}{\mathrm{ch}^\s_{0,n}}
\newcommand{\och}{\ov{\mathrm{ch}}^{\theta^{-1}}_2\!}  

\newcommand{\thetai}{{\theta^{-1}}}

\newcommand{\lN}{\mbox{${\ell}^2(\N_0)$}}

\newcommand{\vare}{\varepsilon}

\newcommand{\ra}{\rightarrow}
\newcommand{\lt}{\triangleleft}
\newcommand{\rt}{\triangleright}
\newcommand{\rti}{\rtimes}
\newcommand{\hs}{\hspace{1pt}}
\newcommand{\hsp}{\hspace{-1pt}}

\newcommand{\tr}{{\mathrm{Tr}}}

\newcommand{\SU}{\cO({\mathrm{SU}}_q(2))}
\newcommand{\PS}{\cO({\mathrm{S}}_q^2)}
\newcommand{\CPS}{\mbox{$\cC^*(\mathrm{S}^2_{qs})$}}
\newcommand{\su}{\U_q({\mathrm{su}}_2)}
\newcommand{\Sq}{\mathrm{S}_q^2}
\newcommand{\BU}{\B\rtimes\U}
\newcommand{\UB}{\U\ltimes\B}
\newcommand{\PSU}{\PS\rtimes\su}

\newcommand{\half}{\mbox{$\frac{1}{2}$}}

\newcommand{\KO}{K_0}

\def\Asnzero{C^\s_0}
\def\Asnone{C^\s_1}
\def\Asntwo{C^\s_2}
\def\Asnthree{C^\s_3}
\def\Asnn{C^\s_n}  \def\ASN{C_\s^n(\A)} \def\ASNN{C_\s^{n-1}(\A)}
\def\Asnnplusone{C^\s_{n+1}}

\def\s{\lambda}
\def\t{\tau}
\def\dd{\mathrm{d}}
\def\ov{\overline}
\def\bo{\boldmath}
\def\mdots{\ldots\hsp}

\begin{document}

\maketitle

\thispagestyle{empty}

\mbox{ }\\[-24pt]
\centerline{{\small Instituto de F\'isica y Matem\'aticas}}
\centerline{{\small Universidad Michoacana de San Nicol\'as de Hidalgo, Morelia, M\'exico }}
%\centerline{{\small 28045 Colima, Colima, M\'exico}}
\centerline{{\small E-mail: Elmar.Wagner@math.uni-leipzig.de}}

\vspace{0.5pc}

\begin{abstract}
The purpose of the paper is twofold:
First, known results of
the noncommutative spin geometry of the standard Podle\'s sphere
are extended by discussing
Poincar\'e duality and orientability. In the discussion of
orientability, Hochschild homology is replaced by a twisted
version which avoids the dimension drop. The twisted Hochschild
cycle representing an orientation is related to the volume form
of the distinguished covariant differential calculus.
Integration over the volume form defines a twisted cyclic 2-cocycle
which computes the $q$-winding numbers of
quantum line bundles.

Second,
a ``twisted'' Chern character from
equivariant $K_0$-theory to even twisted cyclic
homology is introduced which gives rise to a
Chern-Connes pairing between equivariant
$K_0$-theory and twisted cyclic cohomology.
The Chern-Connes pairing between the equivariant
$K_0$-group of the standard Podle\'s sphere
and the generators of twisted cyclic cohomology relative to the
modular automorphism and its inverse are computed.
This includes the pairings with the twisted cyclic 2-cocycle associated
to the volume form,
and the one corresponding to the
``no-dimension drop'' case.
From explicit index computations, it follows that
the pairings with these cocycles
give the $q$-indices of the known equivariant 0-summable Dirac operator
on the standard Podle\'s sphere.
\end{abstract}

\vspace{0.5pc}

\noindent
{\small

\vspace{0.5pc}

}

%\newpage

%*********************************************************
\section{Introduction}                           \label{I}
%*********************************************************

As quantum groups and their associated quantum spaces describe geometric objects
by noncommutative algebras, it is only natural to study them from
Alain Connes' noncommutative geometry point of view \cite{C0}.
The first attempts were made
in the nineties of the past century and exhibited some unexpected features.
For instance, Masuda et al.\ noticed
that the Hochschild dimension of the Podle\'s 2-spheres drops from the
classical dimension 2 to 1 \cite{MNW}.
In another approach, Schm\"udgen proved that
some well-known covariant differential calculi on the quantum SU(2)
cannot be described by a Dirac operator \cite{Sch}.

These observations lowered the expectations on $q$-deformed spaces to be
convincing examples of Connes' noncommutative geometry.
The situation improved after the turn  of the century when the
first spectral triples on $q$-deformed spaces were constructed.
Now there is a lively research activity in studying spectral triples on
quantum groups and their associated quantum spaces.
The best known examples are the (isospectral) spectral triples
on the quantum $\mathrm{SU}(2)$ \cite{CP,DLSSV1} and the 0-dimensional
(i.e., eigenvalues of exponential growth) spectral triple on
the standard Podle\'s sphere \cite{DS} with subsequent analysis of
local index formulas in \cite{C3,DLSSV2} and \cite{NT2}, respectively.
Despite these positive results, some problems remained.
For instance, an equivariant real structure was obtained in
\cite{DLSSV1} only after weakening
the original conditions in \cite{C1}, and the
0-dimensional spectral triple in \cite{DS} is not regular.
To deal with such problems, it was repeatedly suggested to modify the
original axioms of noncommutative spin geometry given in \cite{C2}.

Apart from merely providing examples,  quantum group theory should be
combined with  Connes' noncommutative geometry. The basic input
is equivariance \cite{Sit}, a property which is shared by all the above mentioned
spectral triples. Another substantial step was made by
Kr\"ahmer \cite{K} who constructed Dirac operators on
quantum flag manifolds
and proved that the Dirac operator defines
a finite-dimensional covariant differential calculus in the
sense of Woronowicz \cite{Wo}.
In the case of the standard Podle\'s sphere, considered as $\C P^1$,
Kr\"ahmer's construction reproduces the spectral triple described
by Dabrowski and Sitarz \cite{DS}.
Moreover, Kustermans et al.\ \cite{KMT}
developed a twisted version of cyclic cohomology in order to deal with the
absence of graded traces on quantum groups and Hadfield \cite{H} showed
that the dimension drop can be avoided for the Podle\'s spheres
by considering the twisted Hochschild (co)homology.

Among all the examples mentioned so far, the standard Podle\'s sphere is distinguished.
First of all, because it admits a real structure satisfying
the original conditions in \cite{C1}, and second, the Dirac operator fits nicely
into Woronowicz's theory of  covariant differential calculi \cite{Wo}.
Furthermore, there is a twisted cyclic 2-cocycle associated to the volume
form of the covariant differential calculus \cite{SW3} which reappears in a local
index formula and computes the quantum indices of the Dirac operator \cite{NT2}.
However, from Connes' seven axioms in \cite{C2}, so far only four have been touched.

The main motivation behind the present paper is to expand
the picture of noncommutative spin geometry of the standard Podle\'s
sphere. Throughout the paper, we will work with
the spectral triple found by Dabrowski and Sitarz \cite{DS}.
As indicated above, we partly have to modify the
definitions in \cite{C2}.
The guideline is that the new structures should
still allow the computation of the Chern-Connes and index pairings.

The  Chern-Connes pairing will be defined between twisted cyclic
cohomology and equivariant $\KO$-theory. For this, we recall
the definition of equivariant $\KO$-theory in \cite{NT} and adapt
it to our purpose. This means that we rephrase their definitions
for the left-hand counterpart and take *-str\-uc\-tures into account.
By relating equivariant $\KO$-classes to Hilbert space
representations of crossed product algebras, the equivariant $K_0$-group
of the standard Podle\'s sphere is easily obtained from the results in
\cite{SW2}. More precisely, we show that it is freely generated by
the (equivalence classes of) quantum line bundles of each winding number.

For the  Chern-Connes pairing, we need a ``twisted'' Chern character
mapping equivariant $K_0$-classes into twisted cyclic
homology. Section \ref{TCC} introduces a twisted Chern character
in a general setting. The only requirements are an
appropriate notion of equivariance and a compatibility condition
on the twisting automorphism.

The modular automorphism associated to the Haar state on quantum $\mathrm{SU}(2)$
restricts to an automorphism of the standard Podle\'s sphere, and the
volume form defines a twisted cyclic 2-cocycle relative to it. In Section
\ref{mac}, we compute the full pairing between equivariant $K_0$-theory and twisted
cyclic cohomology with respect to this automorphism.
However, this does not correspond to the ``no-dimension drop'' case of twisted cyclic cohomology.
The dimension drop can be avoided by considering
the inverse modular automorphism. A corresponding
twisted cyclic 2-cocycle, which is also non-trivial on
Hochschild homology, was found by Kr\"ahmer \cite{Kr}.
The Chern-Connes pairing between equivariant $K_0$-theory
and this twisted cyclic 2-cocycle is computed in Section \ref{imac}.
Notably, both twisted cyclic 2-cocycles
compute the $q$-winding numbers.

The discussion of the orientability in Section \ref{O}
demonstrates that, in our example, the twisted versions of
Hochschild and cyclic (co)homology
fit much better into the framework: First, there
is a twisted Hochschild 2-cycle such that its
Hilbert space representation by taking
commutators with the Dirac operator gives the $q$-grading operator,
so the spectral triple satisfies a modified orientability axiom. Second,
the twisted 2-cycle defines a non-trivial class in the twisted
Hochschild homology and also in the twisted
cyclic homology. In particular, it corresponds to the
``no dimension drop'' case. Third, the representation of the
twisted  2-cycle as a 2-form of the algebraically defined
covariant differential calculus
yields the unique (up to a constant) volume form. And finally, integration
over this volume form defines a twisted cyclic 2-cocycle which computes
$q$-indices of the Dirac operator. Note that the combination of the first and
third remark bridges nicely Connes' and Woronowicz' notion of a volume form.
It would be interesting to see whether similar results can be
obtained for other $q$-deformed spaces.

The spectral triple under discussion is not regular. Regularity
provides an operatorial formulation of the calculus of
smooth functions needed in the local index formula.
We do not insist on regularity as long as index computations are
possible. In our example, the indices can be computed by elementary methods
using predominantly equivariance. In Section \ref{IND},
the index and $q$-index of the Dirac operator paired with
any $K_0$-class are calculated.
As expected, we get the winding number (an integer) in the first case
and the $q$-winding number (a $q$-integer) in the second.
Moreover, it is shown that Poincar\'e duality---one of Connes' seven axioms---holds.

Finiteness, another axiom, is implicitly fulfilled by the definition
of the spinor space as projective modules in Section \ref{DO}. We refrain
from the technical details of extending the coordinate algebra to obtain
a pre-C*-algebra.

Although this paper deals only with the standard Podle\'s sphere,
our approach to the Chern-Connes pairing between equivariant
$\KO$-theory and twisted cyclic cohomology
is general enough to be applicable to other examples.
For instance,
Remark \ref{tp} yields a pairing of twisted cyclic cohomology with
equivariant $K_0$-theory, where the only equivariance condition is the
compatibility of the twisting automorphism with the involution of the algebra.
In order not to overstretch the scope of the paper, we did not
consider $K_1$-theory. A definition of a modified $K_1$-group which
uses the modular automorphism can be found in \cite{CPR}.

In this regard, let us remark that our definition of an equivariant $K_0$-group does not
only apply to algebras with a modular automorphism.
Of course, changing the automorphism might change the
equivariant $K_0$-group as much as it might change the
twisted cyclic (co)homology. In the best cases, it will be
isomorphic to the original one, as it happens for the
inverse modular automorphism in the present paper.

%*********************************************************
\section{Preliminaries}                          \label{P}
%*********************************************************

%+++++++++++++++++++++++++++++++++++++++++++++++++++++++++
\subsection{Crossed product algebras}          \label{CPA}
%+++++++++++++++++++++++++++++++++++++++++++++++++++++++++

%
Throughout the paper, we will always work
over the complex numbers $\C$.
Let $\U$ be a Hopf *-al\-ge\-bra
and $\B$ a left $\U$-mo\-dule
*-al\-ge\-bra, that is,  $\B$  is a unital *-al\-ge\-bra with left
$\U$-action $\rt$ satisfying
\begin{equation}                                                 \label{lmod}
f\rt xy = (f_{(1)}\rt x) (f_{(2)}\rt y),\quad
f\rt 1  = \varepsilon(f) 1,\quad
(f\rt x)^* = S(f)^* \rt x^*
\end{equation}
for $x$, $y \in \B$ and $f \in \U$.
Here and throughout the paper, $ \varepsilon$ denotes the counit,
$S$ the antipode, and
$\Delta (f)=f_{(1)} \otimes f_{(2)}$,\,\ $f\in\U$,  is the
Sweedler notation for the comultiplication.

The \emph{left crossed product *-al\-ge\-bra} $\B\rti\U$ is defined as
the *-al\-ge\-bra generated by the two *-sub\-al\-ge\-bras $\B$
and $\U$ with respect to the crossed commutation relations
\begin{equation}                                             \label{fx}
fx = (f_{(1)}\rt x)\hs f_{(2)}, \quad x\in\B,\ \,f\in\U.
\end{equation}

Suppose there exists a faithful state $h$ on $\B$ which
is $\U$-in\-var\-iant, i.e.,
$$
h(f\rt x)=\varepsilon (f)\hs h(x),\quad x\in\B,\ \,f\in\U.
$$
Then there is a unique *-re\-pre\-sen\-ta\-tion $\pi_h$
of  $\BU$ on the domain $\B$ with inner product
$\ip{x}{y}:=h(x^*y)$ such that
\begin{equation}                                           \label{lrep}
\pi_h(y)x=yx,\quad \pi_h(f)x= f\rt x,\quad x,y\in\B,\ \, f\in\U.
\end{equation}

These left-handed definitions have right-handed counterparts.
The right $\U$-action on a right $\U$-mo\-dule *-al\-ge\-bra satisfies
\begin{equation}                                                 \label{rmod}
xy \lt f= (x\lt f_{(1)}) (y\lt f_{(2)}),\quad
1\lt f  = \varepsilon(f) 1,\quad
(x\lt f)^* = x^*\lt  S(f)^*,
\end{equation}
the crossed commutation relations of the
\emph{right crossed product *-al\-ge\-bra} $\UB$ read
\begin{equation}                                          \label{xf}
xf = f_{(1)} \hs (x\lt f_{(2)}) , \quad x\in\B,\ \,f\in\U,
\end{equation}
the invariant state fulfills
$$
h(x\lt f)=\varepsilon (f)\hs h(x),\quad x\in\B,\ \,f\in\U,
$$
and the *-re\-pre\-sen\-ta\-tion $\pi_h$ on $\B$ is given by
\begin{equation}                                           \label{rrep}
\pi_h(y)x=yx,\quad \pi_h(f)x= x\lt S^{-1}(f),\quad x,y\in\B,\ \, f\in\U.
\end{equation}
Proofs of these facts can be found in \cite{SW1} and \cite{SW2}.

%+++++++++++++++++++++++++++++++++++++++++++++++++++++++++
\subsection{Hopf fibration of quantum SU(2)}               \label{HF}
%+++++++++++++++++++++++++++++++++++++++++++++++++++++++++

Throughout this paper, $q$ stands for a positive real number
 such that $q\neq 1$, and we set
$[x]_q:= \frac{q^x-q^{-x}}{q-q^{-1}}$, where $x\in\R$.
For more details on the algebras introduced in this section, we refer
to \cite{KS}.

The Hopf *-al\-ge\-bra $\su$ has four generators
$E$, $F$, $K$, $K^{-1}$ with defining relations
\begin{equation*}                                            %\label{su}
KK^{-1}\hsp=\hsp K^{-1} K\hsp =\hsp 1,\ KE\hsp=\hsp q EK,\
FK\hsp=\hsp qKF,\
EF\hsp-\hsp FE\hsp=\hsp\mbox{$\frac{1}{q-q^{-1}}$} (K^2\hsp-\hsp K^{-2}),
\end{equation*}
involution $E^\ast=F$, $K^\ast=K$, comultiplication
$$\Delta (E)=E\otimes K+K^{-1} \otimes E, \ \
\Delta(F) = F\otimes K+K^{-1} \otimes F,\ \  \Delta (K)=K\otimes K,
$$
counit $\varepsilon (E)\hsp =\hsp\varepsilon (F)\hsp =\hsp
\varepsilon (K-1)\hsp =\hsp 0$,
and antipode  $S(K)\hsp =\hsp K^{-1}$, $S(E)\hsp =\hsp -qE$,
$S(F)\hsp =\hsp -q^{-1}F$.

The coordinate Hopf *-al\-ge\-bra of
quantum $\mathrm{SU}(2)$ will be denoted by $\SU$. A definition of $\SU$ in
terms of generators and relations can be found in \cite{KS}.
Recall from the Peter-Weyl theorem for compact quantum groups that
a linear basis of
$\SU$  is given by the  matrix elements
$t^l_{jk}$ of finite
dimensional unitary corepresentations, where  $l\in\frac{1}{2}\N_0$ and
$j,k=-l,-l+1,\dots,l$.
These matrix elements satisfy
\begin{equation}                                        \label{hopf}
\Delta(t^l_{jk})=\msum_{n=-l}^l t^l_{jn}\otimes t^l_{nk},\quad
\vare(t^l_{jk})= \delta_{jk},\quad S(t^l_{jk}) = t^{l*}_{kj},
\end{equation}
where $\delta_{jk}$ stands for the Kronecker delta.
It follows immediately
that
\begin{equation}                                         \label{t*t}
\msum_{n=-l}^l t^{l*}_{nj}\hs t^l_{nk}=\msum_{n=-l}^l t^{l}_{jn}\hs
t^{l*}_{kn}=\delta_{jk}.
\end{equation}
The standard generators of $\SU$, usually denoted by $a$ and $c$,
are given by $a=t^{1/2}_{-1/2,-1/2}$ and $c=t^{1/2}_{1/2,-1/2}$.
An explicit description of $t^l_{nk}$ in terms of the generators
of $\SU$ can be found in \cite[Section 4.2.4]{KS}.

The Haar state $h$ on $\SU$ is given by $h(t^0_{00})=1$ and
$h(t^l_{jk})=0$ for $l>0$. Since $h$ is faithful, we can define an inner
product on $\SU$ by $\ip{x}{y}:=h(x^*y)$.
With respect to this
inner product, the elements
\begin{equation}                                         \label{v}
v^l_{jk}:=[2l+1]_q^{1/2} q^jt^l_{jk}
\end{equation}
form an
orthonormal vector space basis of $\SU$.

There is a left and a right $\su$-action on $\SU$ turning it into a
$\su$-module *-algebra. Since $h$ is $\su$-invariant,
Equations \eqref{lrep} and \eqref{rrep} define \srep{s} of $\su$.
To distinguish these representations,
we shall omit the representation $\pi_h$ in the first case and write
$\dif_f$ instead of $\pi_h(f)$ in the second case. On the basis vectors
$v^l_{jk}$, the actions of the generators $E$, $F$ and $K$ read
\begin{align}                                            \label{EFK}
E\rt v^{l}_{jk} &= \alpha_k^{l} v^{l}_{j, k+1}, &
F \rt v^{l}_{j k} &= \alpha^{l}_{k-1} v^{l}_{j,k-1}, &
K \rt v^{l}_{jk} &= q^k  v^{l}_{jk},\\                    \label{dif}
\dif_E(v^l_{jk}) &=-\alpha_{j-1}^l v^l_{j-1,k}, &
\dif_F(v^l_{jk}) &=-\alpha_{j}^l v^l_{j+1,k}, &
\dif_K(v^l_{jk}) &=q^{-j} v^l_{jk},
\end{align}
where $\alpha^l_j:=([l-j]_q\hs [l+j+1]_q)^{1/2}$.

Note that, by Equations \eqref{rmod} and \eqref{rrep},
\begin{equation}                                             \label{difab}
\dif_X(ab)= \dif_{X_{(2)}}(a)\dif_{X_{(1)}}(b),\quad
\dif_X(a^*)=\dif_{S(X^{*})}(a)^*
\end{equation}
for $a,b\in\SU$ and $X\in\su$.
We use the right action $\dif_{K^2}$ of the group-like element $K^2$
on $\SU$ to define a Hopf fibration. For $N\in\Z$, set
\begin{equation*}
M_N:=\{ x\in\SU : \dif_{K^2}(x)= q^{-N}x\}
\end{equation*}
and denote by  $\ov M_N$ its Hilbert space closure.
Equation \eqref{dif} implies
\begin{equation*}
M_N=\Lin\{v^l_{N/2,k}: k=-l,\mdots,l,\ \,
l=\mbox{$\frac{|N|}{2}$},\mbox{$\frac{|N|}{2}$}+1,\mdots\}
\end{equation*}
We summarize some basic properties of $M_N$ in the
following lemma.
\begin{lem}                                                   \label{MN}
For $N,K\in\Z$ and $f\in\su$,
\begin{align}                                                 \label{MNMK}
%&\Lin\{a_Nb_K\in M_N M_K\}= M_{N+K},\quad M_N^*\subset M_{-N},\\ \label{EMN}
M_N^*\subset M_{-N},\ \ M_N M_K\subset M_{N+K},&\ \
\Lin\{a_Nb_K\in M_N M_K\}= M_{N+K}, &\\                           \label{EMN}
 \dif_E(M_N)\subset M_{N-2},&\quad \dif_F(M_N)\subset M_{N+2},&\\
 f\rt M_N\subset M_N,&\quad  f\in\su.           &                \label{fMN}
\end{align}
In particular,
$M_0$ is a *-algebra and a left $\su$-module
*-subalgebra of $\SU$, $M_N$ is a $M_0$-bimodule and a left $\su$-module,
and the restriction of the representation $\pi_h$ from \eqref{lrep} to
$M_0$ and $\su$ defines a \srep{} of $M_0\rtimes \su$ on $M_N$.
Moreover, as a left or right $\PS$-module, $M_N$ is generated by
$v^{|N|}_{N,-N},\dots,v^{|N|}_{N,N}$ defined in \eqref{v}.
\end{lem}
\begin{proof}
Most of the assertions are easy consequences of Equations
\eqref{lmod}, \eqref{rmod}, \eqref{EFK} and \eqref{dif}. The last relation
in \eqref{MNMK} follows from a Schur type argument
since $\Lin\{a_Nb_K\in M_N M_K\}\subset M_{N+K}$ is a left $M_0\rtimes \su$-module
and $ M_{N+K}$ is an irreducible one \cite{SW2}. The last claim can be proved by
using an explicit description of the matrix coefficients $t^l_{jk}$ in \eqref{v}
(see, e.g., \cite{KS}).
\end{proof}

The *-algebra $M_0$ is known as the standard
Podle\'s sphere  $\cO(\Sq )$. Usually one defines $\cO(\Sq)$ as the abstract
unital *-al\-ge\-bra with three gen\-e\-ra\-tors $B$, $B^*$,
$A=A^*$ and
defining relations \cite{P}
\begin{equation}                                        \label{algrel}
BA=q^2AB,\ \   AB^* = q^2 B^\ast A,\ \   B^\ast B=A-A^2, \ \
BB^* = q^2 A-q^4 A^2.
\end{equation}
Setting
\begin{equation}                                            \label{ABB*}
A = t^{1/2*}_{1/2,-1/2}\,t^{1/2}_{1/2,-1/2},
\quad B=  t^{1/2*}_{1/2,1/2}\,t^{1/2}_{1/2,-1/2},
\quad B^*=  t^{1/2*}_{1/2,-1/2}\,t^{1/2}_{1/2,1/2},
\end{equation}
yields an embedding into $\SU$ and an isomorphism between $\PS$ and $M_0$
(see, e.g., \cite{KS}).
For generators, the crossed commutation relations \eqref{fx} in $\PSU$
can easily be obtained from \eqref{lmod} and \eqref{EFK}.

Recall that an automorphism $\theta$ satisfying
 $\varphi(xy) = \varphi(\theta (y) x)$
for a state $\varphi$ on a certain *-algebra
is called a modular automorphism (associated to $\varphi$).
It can be shown that $h(xy) = h(\theta (y) x)$ for $x,y\in \SU$, where
\begin{equation}                                            \label{haut}
\theta(y) =\dif_{K^2}(K^{-2} \rt y)= K^{-2} \rt\dif_{K^2}( y).
\end{equation}
The restriction of $h$ to $\PS$ defines a faithful invariant state
on $\PS$ with modular automorphism
\begin{equation}                                            \label{Phaut}
\theta(y) =K^{-2} \rt y,\quad
y\in \PS.
\end{equation}
This follows from \eqref{haut} and the $\dif_{K^2}$-invariance of $\PS$.
By the third relation in \eqref{lmod},
the modular automorphism obeys $\theta(y)^*=\theta^{-1}(y^*)$
for all $y\in\PS$. On the generators
$B$, $B^*$ and $A$, the action of $\theta$  is given by
$$
\theta(B)= q^{2} B, \quad \theta(B^*)= q^{-2} B^*, \quad \theta(A)= A.
$$

%+++++++++++++++++++++++++++++++++++++++++++++++++++++++++
\subsection{Dirac operator on the standard Podle\'s sphere} \label{DO}
%+++++++++++++++++++++++++++++++++++++++++++++++++++++++++

On the standard Podle\'s sphere, there are two non-isomorphic spectral triples known:
the 0-dimensional spectral triple described in \cite{DS} and the isospectral one from
\cite{DDLW}. Both were found by explicit computations on a Hilbert space basis.
However, the 0-dimensional spectral triple admits a convenient description
by using an embedding of the quantum spinor bundle into
the Hopf *-algebra $\SU$ \cite{SW3}.
This construction is unique to the standard Podle\'s sphere;
the general construction of Dirac operators on quantum flag manifolds in \cite{K}
differs slightly from this.

Because of its relation to the representation theory of $\SU$,
we will work in this paper only with the 0-dimensional spectral triple.
The presentation below gives an overview of the results in \cite{SW3} including
simplified ``coordinate free'' proofs.

We define the quantum spinor bundle as the subspace
$$
 W:= M_{-1} \oplus M_{1}\subset \SU
$$
with inner product $\ip{x}{y}=h(x^*y)$, and set
$\hH:= \ov{M}_{-1} \oplus  \ov{M}_{1}$.
The *-re\-pre\-sen\-ta\-tion of $\SU\rtimes \su$ described in \eqref{lrep}
restricts on $W$ to a representation of the crossed product algebra
$\PS\rtimes \su$. For simplicity of notation, we shall omit the symbol $\pi_h$
of the representation.

By Lemma \ref{MN} and Equation \eqref{dif},
the operator
$$
D_0:= \left(
\begin{matrix} 0 & \dif_E \\ \dif_F & 0 \end{matrix} \right)
$$
maps $W$ into itself and  $\{\frac{1}{\sqrt{2}} (v^l_{-1/2 ,k}, \pm  v^l_{1/2 ,k})^\trans:
k=-l,\ldots, l,\ \, l= \frac{1}{2},  \frac{3}{2}, \ldots \}$
forms a complete set of orthonormal eigenvectors.
It follows that the closure of $D_0$ is a
self-adjoint operator, called the Dirac operator $D$. The corresponding eigenvalues
depend only on $l$ and the sign $\pm$, and are given by
$\pm [l+\frac{1}{2}]_q$. In particular, $D$ has compact resolvent.

By \eqref{difab},
$\dif_X(xy)=\dif_{X}(x)\hs \dif_{K^{-1}}(y)+\dif_{K}(x)\hs \dif_{X}(y)$
for $x,y\in \SU$ and $X=E,F$.
Using  $\dif_{K}(v_N)= q^{-N/2}v_N$ for $v_N\in M_N$, one computes
\begin{equation}                                                     \label{Da}
 [D,a]=\left(\begin{matrix} 0 & q^{1/2}\dif_E(a) \\ q^{-1/2}\dif_F(a) & 0
 \end{matrix} \right),\quad  a\in\PS.
\end{equation}
As a consequence, $ [D,a]\in \mathrm{B}(\hH)$ for all $a\in\PS$.

There is a natural grading operator $\gamma$ on $\hH$ given by
$$
\gamma:= \left(
\begin{matrix} 1 & 0 \\ 0 & -1 \end{matrix} \right).
$$
Clearly, $D\gamma=-\gamma D$ and $\gamma a=a\gamma$ for all
$a\in\PS$.

Set $J_0:=*\circ\dif_K K^{-1}$.
For $x,y\in\SU$,
\begin{align*}
\ip{J_0(x)}{J_0(y)}\!&=\!h\big((\dif_K K^{-1}\rt x)(\dif_{K^{-1}} K\rt y^*)\big)
\!=\!h\big((\dif_{K} K^{-1}\rt y^*)(\dif_K K^{-1}\rt x)\big)\\
\!&=\!h\big(\dif_{K} K^{-1}\rt (y^*x)\big)    \!=\!h(y^*x)\!=\!\ip{y}{x}
\end{align*}
by Equations \eqref{lmod}, \eqref{rmod}, \eqref{haut}, and the
$\su$-invariance of $h$. Hence $J_0$ is anti-unitary.
From  $*\circ\dif_K K^{-1}=\dif_{K^{-1}} K \circ *$, we conclude
$J_0^2=1$ and
\begin{equation}                                      \label{J0}
J_0\hs a^*J_0^{-1}x=J_0\hs a^*J_0x= \dif_{K^{-1}} K\rt\big(a^* (\dif_K K^{-1}\rt x)^*\big)^*
=x (K\rt a)
\end{equation}
for $a\in\PS$ and $x\in\SU$. By Lemma \ref{MN}, $J_0$ leaves $W$ invariant.
Since the representations of $\PS$ and $[D,a]$, $a\in\PS$, act on $W$ by
left multiplication, it follows from \eqref{J0} that
$$
[a,J_0\hs b^*J_0^{-1}]=0,\quad \big[[D,a],J_0\hs b^*J_0^{-1}\big]=0,\quad a,b\in\PS.
$$

The last relations in \eqref{rmod} and \eqref{rrep} imply
$J_0\hs \dif_f= \dif_{S(K^{-1} f^* K)} J_0$ and therefore
$J_0 D=-DJ_0$ on $W$. Clearly, $\gamma J_0=-J_0\gamma$ on $W$ by \eqref{MNMK}.
Considering $J:=\gamma J_0$ as an anti-unitary operator on $\hH$, one gets $J^2=-1$ and
$JD=DJ$. Summing up, we obtain the following theorem \cite[Theorem 3.3 (iii)]{SW3}.

\begin{thm}
The quintuple $(\PS,\hH,D,J,\gamma)$ is a real even spectral triple
in the sense of \cite{C1}.
\end{thm}

The formulas $\Omega:=\Lin\{b\,[D,a] : a,b\in\PS\}$ and
$\dd a:=\im [D,a]$ define a covariant first order differential
calculus $\dd:\PS \ra \Omega$. The corresponding universal differential
calculus is given by
$\Omega^{\wedge} = \oplus_{k=0}^\infty \Omega^{\otimes k} / \cI$, where
$\Omega^{\otimes k}= \Omega\otimes_{\PS}\cdots \otimes_{\PS}\Omega$
($k$-times),  $\Omega^{\otimes 0}=\PS$,
and $\cI$ denotes the two-sided ideal in the tensor algebra
$\oplus_{k=0}^\infty \Omega^{\otimes k}$ generated by the elements
$\sum_i \dd a_i \otimes \dd b_i$ such that $\sum_i  a_i\, \dd b_i=0$.
The product in the algebra $\Omega^{\wedge}$ is denoted by $\wedge$ and
we write
$\dd(a_0\dd a_1\wedge\cdots\wedge\dd a_k)
=\dd a_0\wedge\dd a_1\wedge\cdots\wedge\dd a_k$.
The following facts are
proved in \cite{SW3}.
\begin{prop}                                                \label{cycle}
There exists an invariant 2-form $\omega\in\Omega^{\wedge 2}$ such that
$\Omega^{\wedge 2}$ is the free $\PS$-module generated by $\omega$
with $a\omega=\omega a$ for all $a\in\PS$ and
$a_0\dd a_1\wedge\dd a_2
=a_0\big(
q\dif_E(a_1)\dif_F(a_2)
-q^{-1}\dif_F(a_1)\dif_E(a_2)
\big)\omega$.
Let $h$ denote the Haar state on $\PS$ and $\theta$ its modular
automorphism described in \eqref{Phaut}.
The linear functional $\int : \Omega^{\wedge 2}\ra \C$,\ \,$\int a\omega := h(a)$
defines a non-trivial  $\theta$-twisted cyclic 2-cocycle $\tau$ on $\PS$ given by
\begin{equation}                                             \label{cyc}
\tau(a_0,a_1,a_2)=\int a_0\,\dd a_1\wedge\dd a_2
=h\Big(a_0\big(q\hs\dif_E(a_1)\dif_F(a_2)-q^{-1}\dif_F(a_1)\dif_E(a_2)\big)\Big).
\end{equation}
\end{prop}
We call $\omega$ a volume form associated with the covariant differential
calculus. For a definition of twisted cyclic cocycles, see Section \ref{TCC}.
An explicit expressions of $\omega$ can easily be deduced from the
formulas given in Section \ref{O} and in \cite[Appendix (Proof of Lemma 4.4)]{SW3}.

%*********************************************************
\section{Equivariant \bo{$\KO$}-theory}                  \label{EKT}
%*********************************************************

\subsection{Definition and basic material}                \label{DB}

This section is concerned with a simple definition of
equivariant $\KO$-theory. The approach follows closely the lines of
\cite{NT} which works well for compact quantum groups.
Our treatment differs from that in \cite{NT} in two aspects.
First, we take *-structures into account, and second, we will include also
{\it left} crossed product algebras in our considerations.

We start by recalling some definitions from \cite{NT}.
Let $\B$ be a right $\U$-module algebra. Suppose that
 $\rho^\circ:\U^\circ\ra \End(\C^n)$ is a finite dimensional
 representation of the opposite algebra $\U^\circ$ or, equivalently,
 $\rho^\circ:\U\ra \End(\C^n)$ is a finite dimensional
 \emph{anti}-homomorphism.
 Then $\C^n\otimes\B$ inherits a right $\U\ltimes\B$-module structure from the
 anti-representation
 \begin{equation}                                             \label{Rrep}
 \pi^\circ(f)(v\otimes a) :=  \rho^\circ(f_{(1)})v \otimes a\lt f_{(2)} ,\quad
\pi^\circ(b)(v\otimes a) := v\otimes a b,\quad f\in\U, \,b\in\B.
\end{equation}
The algebra
$\M_{n\times n}(\C)\otimes \B$
can be embedded into $\End(\C^n\otimes \B)$ by
$$
   \M_{n\times n}(\C)\otimes\B  \,\ni\, T\otimes b\longmapsto
   T\otimes \mL_b\,\in\,\End(\C^n\otimes \B),
$$
where $\mL_b\hs a:= ba$, $\,a,b\in\B$,
denotes the left multiplication of $\B$.
On column vectors $\fett{b}\in \B^n \cong \C^n\otimes \B$,  the action of
$X\in  \M_{n\times n}(\B) \cong  \M_{n\times n}(\C)\otimes \B$
is conveniently expressed by matrix multiplication, i.e.,
\begin{equation}                                                \label{Xb}
   \M_{n\times n}(\B)\,\ni\, X \longmapsto
\big( \fett{ b} \mapsto X\,\fett{ b}\big) \,\in\,\End(\C^n\otimes\B ).
\end{equation}
We turn $\End(\C^n\otimes\B )$
into a right $\U$-module by using the left adjoint action
of $\U^\circ$, i.e.
\begin{equation*}
 \ad_{\mathrm{L}}^\circ(f)\hs (X):= \pi^\circ(f_{(1)})\hs X\hs \pi^\circ(S^{-1}(f_{(2)})),
\quad   X\in\End(\C^n\otimes\B ),\ \, f\in\U.
\end{equation*}
From \cite[Lemma 1.1]{NT} (or from direct calculations), it follows that
\begin{equation}                                              \label{adoL}
\ad_{\mathrm{L}}^\circ(f)\hs (T \otimes \mL_b) =
\rho^\circ(f_{(1)})T\rho^\circ(S^{-1}(f_{(3)}))
\otimes \mL_{b\lt f_{(2)}}
\end{equation}
for all $T\otimes b\in   \M_{n\times n}(\C)\otimes\B$ and $f\in\U$.
Under the identification \eqref{Xb}, Equation \eqref{adoL} becomes
\begin{equation}                                                \label{Xf}
  \ad_{\mathrm{L}}^\circ(f)(X)=
  \rho^\circ(f_{(1)})\big(X\lt f_{(2)}\big)\rho^\circ(S^{-1}(f_{(3)})),
\end{equation}
where $X\in  \M_{n\times n}(\B)$,\ \,$f\in\U$, and
$X\lt f$ stands for the action of $f$ on each entry of the matrix $X$.
Looking at the last equations, one readily sees that $\M_{n\times n}(\B) $ is
a right $\U$-module subalgebra of $\End(\C^n\otimes\B)$.
Alternatively, we can consider $\M_{n\times n}(\B) $ as
a left  $\U^\circ$-module subalgebra of $\End(\C^n\otimes\B)$.

Now we take *-structures into account.
Suppose that
$\U$ is a Hopf *-algebra, $\B$ a right $\U$-module *-algebra,
and $\rho^\circ :\U^\circ\ra\End(\C^n)$ a \srep.
In order to turn $\M_{n\times n}(\B)$ into a $\U^\circ$-module *-algebra,
we need an automorphism $\sigma:\B \ra \B$ such that
$\sigma(a\lt f )=\sigma(a)\lt S^{-2}(f)$ and $\sigma(a)^*=\sigma^{-1}(a^*)$
for all $a\in \B$ and $f\in \U$.
Then $X^\dagger:= \sigma(X)^*$ defines an involution on
$\M_{n\times n}(\B)$ such that
\begin{align*}
 \big(\ad_{\mathrm{L}}^\circ(f)(X)\big)^\dagger &=
 \rho^\circ(S^{-1}(f_{(3)})^*)\big(\sigma(X)^*\lt S^{-1}(f_{(2)})^*\big)
 \rho^\circ(f_{(1)}^*)\\
 &=\ad_{\mathrm{L}}^\circ(S^{-1}(f)^*)(X^\dagger)
\end{align*}
for all $X\in \M_{n\times n}(\B)$ and $f\in\U$.
As $S^{-1}$ is the antipode of $\U^\circ$, the last relation shows that
$\M_{n\times n}(\B)$ with the involution ${}^\dagger$ and
the left $\U^\circ$-action $\ad_{\mathrm{L}}^\circ$ is
a left $\U^\circ$-module *-algebra.

A matrix $X\in  \M_{n\times n}(\B)$ is called
\emph{(right) $\U$-invariant} if there exists a \srep{} $\rho^\circ :\U^\circ\ra\End(\C^n)$
such that
$$
\ad_{\mathrm{L}}^\circ(f)(X)=\vare(f) X
$$
for all $f\in\U$.

Next we make analogous definitions for left crossed product algebras.
Let thus $\B$ be a {\it left} $\U$-module algebra.
In order to apply the definitions given above, we consider $\B$ as
a {\it right} $\U^\cop$-module with right $\U^\cop$-action given by
$a\lt f:=S^{-1}(f)\rt a$, where $f\in \U$ and $a\in\B$.
Then $\ad_{\mathrm{L}}^\circ$ turns $\End(\C^n\otimes\B )$
into a left  $\U^{\circ,\cop}$-module algebra.
To get back to a $\U^{\circ}$-module algebra, we use again the
inverse of the antipode and define a right $\U^{\circ}$-action
on $\End(\C^n\otimes\B )$ by setting
$\ad_{\mathrm{R}}^\circ(f):=\ad_{\mathrm{L}}^\circ(S^{-1}(f))$.
Now $\M_{n\times n}(\B)$ becomes a right $\U^{\circ}$-module
subalgebra and, for all $X\in \M_{n\times n}(\B)$,
\begin{equation}                                                \label{XfL}
  \ad_{\mathrm{R}}^\circ(f)(X)=
  \rho^\circ(S^{-1}(f_{(1)}))\big(S^{-2}(f_{(2)})\rt X\big)\rho^\circ(f_{(3)}).
\end{equation}
Similarly to the above, we assume that there is an automorphism
$\sigma:\B \ra \B$ such that
$\sigma(f\rt a )=S^{2}(f)\rt \sigma(a)$ and $\sigma(a)^*=\sigma^{-1}(a^*)$
for all $a\in \B$ and $f\in \U$. With respect to the involution
$X^\dagger:= \sigma(X)^*$, we get
\begin{align*}
 \big(\ad_{\mathrm{R}}^\circ(f)(X)\big)^\dagger &=
 \rho^\circ(f_{(3)}^*)\big( S(f_{(2)})^*\rt \sigma(X)^*\big)
 \rho^\circ(S^{-1}(f_{(1)})^*)\\
 &=\ad_{\mathrm{R}}^\circ(S^{-1}(f)^*)(X^\dagger)
\end{align*}
for all $X\in \M_{n\times n}(\B)$ and $f\in\U$
since $S^{-1}\circ *=*\circ S$. Hence the involution~${}^\dagger$ and the
right $\U^\circ$-action $\ad_{\mathrm{R}}^\circ$ endow $\M_{n\times n}(\B)$
with the structure of a right $\U^\circ$-module *-algebra.
Note that the automorphism $S^{-2}$ in Equation \eqref{XfL} is necessary
for $(\ad_{\mathrm{R}}^\circ(f)(X))^\dagger
=\ad_{\mathrm{R}}^\circ(S^{-1}(f)^*)(X^\dagger)$ to hold.

As above, we say that $X\in  \M_{n\times n}(\B)$ is
\emph{(left) $\U$-invariant} if there exists a \srep{} $\rho^\circ :\U^\circ\ra\End(\C^n)$
such that
$$
\ad_{\mathrm{R}}^\circ(f)(X)=\vare(f) X
$$
for all $f\in\U$.

For a definition of equivariant $\KO$-theory, we shall use the
Murray-von Neumann equivalence of projections.
Given an automorphism $\sigma$ of $\B$
such that  $\sigma(b)^*=\sigma^{-1}(b^*)$,
an idempotent $P\in \M_{n\times n}(\B)$ will be called
{\it projection} if $P=P^\dagger$.

\begin{defn}                                                 \label{MvN}
Let $\B$ be a *-algebra and $\sigma:\B \ra \B$ an automorphism
satisfying $\sigma(b)^*=\sigma^{-1}(b^*)$.
Suppose that $\B$ is a right $\U$-module *-algebra and
$\sigma(a\lt f )=\sigma(a)\lt S^{-2}(f)$ for all
$a\in \B$ and $f\in \U$; or
$\B$ is a left $\U$-module *-algebra and
$\sigma(f\rt a )=S^{2}(f)\rt \sigma(a)$.

For $n,m\in\N$, let $\rho_1^\circ: \U^\circ \ra \End(\C^n)$ and
$\rho_2^\circ: \U^\circ \ra \End(\C^m)$ be finite dimensional
*-re\-pre\-sen\-ta\-tions.
Denote by $\pi_1^\circ$ and $\pi_2^\circ$ the representations of
$\U^\circ$ on $\C^n\otimes \B$ and $\C^m\otimes \B$, respectively, given
in Equation \eqref{Rrep} or \eqref{Lrep}. We say that invariant projections
$P\in\M_{n\times n}(\B)$  and $Q\in\M_{m\times m}(\B)$
are \emph{Murray-von Neumann equivalent} if there exists a
$V\in \Hom_\B( \C^n\otimes\B,\C^m\otimes \B)$
such that $V^\dagger V= P$, $VV^\dagger=Q$ and
$\pi_2(f)V=V\pi_1(f)$ for all
$f\in\U$.
\end{defn}
\begin{rem}
Since $\Hom_\B(\C^n\otimes\B ,\C^m\otimes\B)\cong
\M_{m\times n}(\B)$,
we can assume that $V\in\M_{m\times n}(\B)$.
\end{rem}
We are now in a position to state the following practical definition
of equivariant $\KO$-theory.
\begin{defn}                                                \label{K}
Let $\U$, $\B$ and $\sigma$ be as in Definition \ref{MvN}.
Then $\KO^{\U}(\B)$  (resp.\ ${}^{\U}\hsp\KO(\B)$)
denotes the Grothendieck
group obtained from the additive semigroup
of Murray-von Neumann equivalent,
$\ad_{\mathrm{R}}^\circ$ (resp.\ $\ad_{\mathrm{L}}^\circ$) invariant
 projections of any size of $\B$-valued square matrices with addition
$[P]+[Q] = [P\oplus Q]$ and additive identity $0=[0]$
for any size of zero matrix.
\end{defn}
\begin{rem}                                               \label{os}
If $\rho_1^\circ: \U^\circ\ra \End(\C^n)$ and
$\rho_2^\circ: \U^\circ\ra \End(\C^m)$ are finite dimensional
*-re\-pre\-sen\-ta\-tions,
and $P\in\M_{n\times n}(\B)$  and $Q\in\M_{m\times m}(\B)$
are invariant projections, then the notation
$P\oplus Q$ refers to the projection $\mathrm{diag}(P,Q)$
in $\M_{(n+m)\times (n+m)}(\B)\cong \End(\C^n\oplus \C^m) \otimes\B$,
where
the representation of $\U^\circ$ on $\C^n\oplus \C^m$ is given by
$\rho_1^\circ\oplus \rho_2^\circ$.
\end{rem}
\begin{ex}                                            \label{Ks}
Suppose that $\B$ is a unital *-algebra and $\sigma$ an
automorphism satisfying $\sigma(a)^*=\sigma^{-1}(a^*)$.
Then we can define a commutative and co-commutative
Hopf *-algebra $\U(\sigma)$ generated by $\sigma$
with Hopf structure
$$
\Delta(\sigma)=\sigma\otimes \sigma,\quad \vare(\sigma)=1,\quad
S(\sigma)=\sigma^{-1}
$$
and involution $\sigma^*=\sigma$. The left and right actions
$$
\sigma\rt b= b\lt \sigma =\sigma(b),\qquad b\in\B,
$$
turn $\B$ into a left and right $\U(\sigma)$-module *-algebra such that,
for all $f\in \U$, $\sigma(b\lt f )=\sigma(b)\lt S^{-2}(f)$ and
$\sigma(f\rt b )=S^{2}(f)\rt \sigma(b)$ since $S^{2}=\id$. In this way we
obtain a definition of equivariant $\KO$-theory which
depends only on the automorphism $\sigma$. Instead of
$\KO^{\U(\sigma)}\hsp (\B)\ \,\big(={}^{\U(\sigma)}\hsp\KO(\B)\big)$,
we shall from now on simply
write $\KO^\sigma\hsp(\B)$.

This definition of equivariant $\KO$-theory is
strongly related to $\sigma$-twisted cyclic (co)homology.
In particular, as we shall see in Remark \ref{tp},
it allows us to define a pairing between
$\KO^\sigma(\B)$ and twisted cyclic cohomology.
\end{ex}

%++++++++++++++++++++++++++++++++++++++++++++++++++++++
\subsection{Equivariant \bo{$\KO$}-theory and
the modular automorphism}                                \label{sec-equi-mod}
%++++++++++++++++++++++++++++++++++++++++++++++++++++++

In this section we show that, in presence of a modular
modular automorphism, equivariant $\KO$-classes  are
intimately related to unitarily equivalent
Hilbert space representations of the opposite
crossed product algebra. The computation of the $\KO$-group
of the standard Podle\'s sphere
can then be reduced to the classification of
certain types of unitarily equivalent Hilbert space representations.
The details below give also an a posteriori motivation for
the definitions made in the previous section.

Throughout this section, we suppose that
$\B$ is a left (or right) $\U$-module *-algebra and
$h:\B\ra\C$ is a faithful invariant state
with modular automorphism $\theta$.
Recall from Section \ref{CPA} that $\ip{a}{b}:=h(a^*b)$ defines an inner product
on $\B$ such that
$h(a^*(f\rt b))=\ip{a}{f\rt b} =\ip{f^*\rt a}{b} =h((f^*\rt a)^*b)$
for all $a,b\in\B$ and $f\in\U$.
From
$$
h(\theta(a^*)b)=h(ba^*)=\overline{h(ab^*)}=\overline{h(b^*\theta^{-1}(a))}
=h(\theta^{-1}(a)^*b),
$$
it follows that $\theta(a^*)=\theta^{-1}(a)^*$, and
\begin{align*}
h(\theta(f\rt b)\hs a^*)&= h(a^*\hs (f\rt b))
= h((f^*\rt a)^* \hs b) = h(\theta(b)\hs  (S(f^*)^*\rt a^*))\\
&= h(( S(f^*)\rt\theta(b)^*)^* \hs  a^*)
= h(( S^{-2}(f)\rt\theta(b)) \hs  a^*),
\end{align*}
implies $\theta(f\rt b))=S^{-2}(f)\rt\theta(b)$.
Similarly one shows that $\theta(b\lt f))=\theta(b)\lt S^{2}(f)$.
Hence $\sigma:=\theta^{-1}$ satisfies the conditions of Definition \ref{MvN}.

Next we define an inner product on $\C^n\otimes\B$ by
\begin{equation}                                       \label{ip}
\ip{v\otimes a}{w\otimes b}^\circ := h(b\hs a^*)\ip{v}{w}_{\hsp_{\C^n}}.
\end{equation}
Note that we used $h(ba^*)$ instead of $h(a^*b)$ so that the
right multiplication yields a *-representation of the
opposite algebra $\B^\circ$. If $\B$ is a right $\U$-module *-algebra
and $\rho^\circ :\U^\circ\ra\End(\C^n)$ is a \srep{}, then
Equation \eqref{Rrep} defines a \srep{} of the
opposite crossed product algebra $(\UB)^\circ$.

Similarly, if $\B$ is a left $\U$-module *-algebra, we set
\begin{equation}                                             \label{Lrep}
 \pi^\circ(f)(v\otimes a) :=  \rho^\circ(f_{(2)})v \otimes S^{-1}(f_{(1)})\rt a  ,\quad
\pi^\circ(b)(v\otimes a) := v\otimes a b,
\end{equation}
where $b\in\B$,  $f\in\U$ and $v\otimes a\in\C^n\otimes\B$.
One easily checks that Equation \eqref{Lrep} defines a \srep{}
of opposite crossed product algebra $(\BU)^\circ$, where the
inner product on $\C^n\otimes\B$ is given by \eqref{ip}.

By Equation \eqref{Xb}, matrix multiplication from the left defines an
embedding of $\M_{n\times n}(\B)$ into $\End(\C^n\otimes \B)$.
Given $X=(x_{ij})_{i,j=1}^n\in \M_{n\times n}(\B)$, let $X^+$ denote the Hilbert
space adjoint of the corresponding operator in $\End(\C^n\otimes \B)$.
Then, for all $\fett{ a}=(a_1,\mdots, a_n)^\trans\in \B^n$
and $\fett{ b}=(b_1,\mdots, b_n)^\trans\in \B^n$,
\begin{align*}
\ip{ \fett{ a}}{X^+\fett{ b}}^\circ &= \ip{X\hs\fett{ a}}{ \fett{ b}}^\circ
= \msum_{i,j=1}^n h(b_j a^*_i x_{ji}^* )
= \msum_{i,j=1}^n h(\theta(x_{ji}^*)b_ja^*_i) )\\
&=\ip{ \fett{ a}}{\theta(X^*)\fett{ b}}^\circ .
\end{align*}
Thus $X^+=\theta(X^*)=\sigma(X)^*=X^\dagger$ and the above embedding
becomes a \srep. In particular, projections in $\M_{n\times n}(\B)$ yield
orthogonal projections on $\C^n\otimes \B$.

The following proposition relates invariant projections
to Hilbert space representations
of the opposite crossed product algebra.

\begin{prop}                                            \label{idem}
Let $P\in\M_{n\times n}(\B)$ be a projection.
Then the restriction of $\pi^\circ$ to the
projective right $\B$-module $P\hs\B^n$ defines a \srep\
of $(\UB)^\circ$ (or $(\BU)^\circ$) if and only if $P$ is invariant.
\end{prop}
\begin{proof}
We prove Proposition \ref{idem} for left $\U$-module *-algebras, the
proof for the right-handed counterpart is similar.

For column vectors $\fett{ b}=(b_1,\mdots, b_n)^\trans\in \B^n\cong
\C^n\otimes\B $,
the first equation of \eqref{Lrep} can be written
$\pi^\circ(f)\fett{ b}= \rho^\circ(f_{(2)})\hs(S^{-1}(f_{(1)}) \rt\fett{ b} )$.
Thus
\begin{align*}
\pi^\circ(f)( P\hs \fett{ b}) &= \rho^\circ(f_{(3)})\hs (S^{-1}(f_{(2)})\rt P)\hs
(S^{-1}(f_{(1)})\rt \fett{ b})\\
&=\rho^\circ(f_{(5)})\hs (S^{-1}(f_{(4)})\rt P)\hs\rho^\circ(S(f_{(3)}))\hs \rho^\circ(f_{(2)})
\hs(S^{-1}(f_{(1)})\rt\fett{ b} ) \\
&= \ad_{\mathrm{R}}^\circ(S(f_{(2)}))(P)\, \pi^\circ(f_{(1)})\hs \fett{ b}.
\end{align*}
Hence, as Hilbert space operators on $\B^n$,
$$
\pi^\circ(f)\hs P = \ad_{\mathrm{R}}^\circ(S(f_{(2)}))(P)\hs\hs \pi^\circ(f_{(1)}).
$$
Since $\ad_{\mathrm{R}}^\circ(S(f))(P)=\ad_{\mathrm{L}}^\circ(f)(P)
=\pi^\circ(f_{(1)})\hs P\hs \pi^\circ(S^{-1}(f_{(2)}))$, it follows that
$\pi^\circ(f) P = P\hs \pi^\circ(f)$ if and only if
$\ad_{\mathrm{R}}^\circ(f)(P) = \vare(f) P$.
\end{proof}
%

%++++++++++++++++++++++++++++++++++++++++++++++++++++++
\subsection{The equivariant \bo{$\KO$}-group of
the standard Podle\'s sphere}                                \label{sec-equiv}
%++++++++++++++++++++++++++++++++++++++++++++++++++++++

We restrict ourselves to the computation of $\KO^{\su}\hsp(\PS)$
with respect to the inverse modular automorphism $\theta^{-1}$ of $\PS$,
the computation of its right-handed counterpart is analogous.
In particular, the outcome would be that
$\KO^{\su}\hsp(\PS)\cong{}^{\su}\hsp\KO(\PS)$.

Our first aim is to construct representatives for
equivariant $\KO$-classes.
For $n\in\half\Z$ and $l=|n|,|n|+1,\ldots\,$, let $\fett{t}^l_n$
denote the row vector
\begin{equation}                                          \label{rv}
\fett{t}^l_n:=(t^l_{n,-l},t^l_{n,-l+1},\mdots,t^l_{n,l}),
\end{equation}
where $t^l_{n,k}$ are the matrix elements from Section \ref{HF}.
Recall from Section \ref{HF} that there is a  $(2l+1)$-dimensional
\srep{} of $\su$ on
$$
V^l_n:=\Lin\{t^l_{n,k}:k=-l,\mdots,l\}=\Lin\{v^l_{n,k}:k=-l,\mdots,l\},\quad l\geq |n|,
$$
given  by Equation \eqref{EFK}.
These representations are irreducible and called spin-$l$-representations.
Let $\sigma_l:\su \ra \M_{(2l+1)\times (2l+1)}(\C)$ be the matrix
representation determined by
$f\rt \msum_{k=-l}^l \alpha_k v^l_{nk}
=\msum_{k,j=-l}^l \sigma_l(f)_{jk} \alpha_k v^l_{nj}$.
Then
\begin{equation}                                            \label{ft}
f\rt \fett{t}^l_n= \fett{t}^l_n\hs \sigma_l(f),\quad
f\rt \fett{t}^{l*}_n=\sigma_l(S(f))\hs\fett{t}^{l*}_n, \quad f\in\su,
\end{equation}
where we used \eqref{lmod} in the second relation.
Consider the homomorphism
\begin{equation}                                            \label{rhoo}
\rho_l^\circ: \su^\circ \ra \End(\C^{2l+1}),\quad
\rho_l^\circ(f):= \sigma_l(K^{-1}S(f)K).
\end{equation}
From $K^2 f K^{-2}=S^2(f)$ and $S(f)^*= S^{-1}(f^*)$ for all $f\in\su$,
it follows that $\rho_l^\circ$ is a \srep{} of $\su^\circ$.

Finally, for $N\in\Z$, define
\begin{equation}                                                   \label{PN}
 P_{N}:= \rho_{|N|/2}^\circ(K^{-1})\, \fett{t}^{|N|/2*}_{N/2}\,\fett{t}^{|N|/2}_{N/2}
\,\rho_{|N|/2}^\circ(K).
\end{equation}
We summarize some
crucial properties of these matrices  in the next lemma.
\begin{lem}                                                       \label{PLN}
The matrices $P_N$ belong to $\M_{|N|+1\times |N|+1}(\PS)$
and are $\adro$-invariant projections with respect to
the anti-representation $\rho_{|N|/2}^\circ$ of $\su$ and
the involution  $P_N^\dagger=\theta^{-1}( {P_N})^*$, where
$\theta^{-1}(b)=K^{2}\rt b$ for all $b\in\PS$.
\end{lem}
\begin{proof}
For brevity of notation, set $n:=N/2$ and $l:=|N|/2$.
From
$$
\dif_{K^2} (t^{l*}_{n,j} t^l_{n,k})
\,=\,\big(\dif_{K^{-2}} (t^{l}_{n,j})\big)^* \dif_{K^2} (t^l_{n,k})
\,=\, t^{l*}_{n,j} t^l_{n,k},
$$
it follows that the entries of  $P_N$ belong to $\PS$.
Using  $K^2 f K^{-2}=S^2(f)$ and Equations \eqref{XfL} and \eqref{ft}--\eqref{PN}, we get
\begin{align*}
 \adro(f)&(P_N)= \rho_l^\circ(S^{-1}(f_{(1)}))
\big( S^{-2}(f_{(2)})\rt P_N\big)\rho_l^\circ(f_{(3)})\\
&=\sigma_l(K^{-1}f_{(1)}K^2) (S^{-2}(f_{(2)})\rt\fett{t}^{l*}_{n})
(S^{-2}(f_{(3)})\rt\fett{t}^{l}_{n}) \sigma_l(K^{-2}S(f_{(4)})K)\\
&=\sigma_l(K^{-1}f_{(1)}K^2 S^{-1}(f_{(2)}))\,\fett{t}^{l*}_{n}\,
\fett{t}^{l}_{n}\, \sigma_l(S^{-2}(f_{(3)})K^{-2}S(f_{(4)})K)\\
&=\sigma_l(K^{-1}f_{(1)} S(f_{(2)})K^2)\,\fett{t}^{l*}_{n}\,
\fett{t}^{l}_{n}\, \sigma_l(K^{-2}f_{(3)}S(f_{(4)})K)\\
&=\vare(f)\,\sigma_l(K)\,\fett{t}^{l*}_{n}\,
\fett{t}^{l}_{n}\, \sigma_l(K^{-1}) = \vare(f)\, P_N
\end{align*}
for all $f\in\su$. Thus $P_N$ is $\adro$-invariant.
The $\adro$-invariance implies
\begin{equation*}
\theta^{-1}(P_N)=K^2\rt P_N
= \rho_l^\circ(K^2)\,\adro(K^2)(P_N) \,\rho_l^\circ(K^{-2})
= \rho_l^\circ(K^2)\, P_N \,\rho_l^\circ(K^{-2}),
\end{equation*}
hence $P_N^\dagger := \theta^{-1}(P_N)^*
= \big(\rho_l^\circ(K)\, \fett{t}^{l*}_{n}\,\fett{t}^l_{n} \,
\rho_l^\circ(K^{-1}) \big)^*=P_N$.
Clearly, $P_N^2=P_N$ since $\fett{t}^{l}_{n}\,\fett{t}^{l*}_{n}=1$
by \eqref{hopf},
so $P_N$ is a projection.
\end{proof}

By Definition \ref{K}, the projections in the last
lemma represent classes of $\KO^{\su}\hsp(\PS)$
relative to the inverse modular automorphism $\theta^{-1}$ of $\PS$.
The next proposition shows that
the projections $P_N$ determine completely
$\KO^{\su}\hsp(\PS)$.
\begin{prop}                                                       \label{eqK0}
With respect to the automorphism $\theta^{-1}:\B\ra\B$ given by
$\theta^{-1}(b)=K^2\rt b$,
the group $\KO^{\su}\hsp(\PS)$ is isomorphic to the free abelian
group with one generator, $[P_N]$, for each $N\in\Z$.
\end{prop}
\begin{proof}
By Proposition \ref{idem},
the projections $P_N$ determine Hilbert space representations of
$(\PSU)^\circ$. Our first aim is to show that, for $N\in\Z$,
we obtain pairwise inequivalent integrable representations and that
each Hilbert space representations arising from
$\adro$-invariant $\PS$-valued projections
decomposes into those obtained by $P_N$.
Here, a \srep{} of $(\PSU)^\circ$ is called
\emph{integrable}
if its restriction to $\su^{\circ}$ decomposes into finite dimensional
\srep{}s. Note that, by the Clebsch-Gordon decomposition,
the (tensor product) representation of $(\PSU)^\circ$ on
$\C^n\otimes\PS $ defined in \eqref{Lrep} is
integrable and so are the representations from Proposition \ref{idem}.

Straightforward calculations show that $(\PSU)^\circ$
is isomorphic to
$\cO({\mathrm{S}}_{q^{-1}}^2)\rtimes\U_{q^{-1}}({\mathrm{su}}_2)$
with $A$ replaced by $q^{-2}A$.
By using alternatively the opposite algebras and the replacement
$q\mapsto q^{-1}$, we can apply freely the results from \cite{SW2}.

To begin, consider the Hopf fibration of $\SU^\circ$ given by
\begin{equation*}
\SU^\circ=\oplus_{N\in\Z}M_N^\circ,\quad
M_N^\circ:=\{ x\in\SU^\circ : \dif_{K^2}(x)= q^{-N}x\}.
\end{equation*}
As in \cite{SW2}, the restriction of the GNS-re\-pre\-sen\-ta\-tion $\pi_h$
to $M_N^\circ$ defines pairwise inequivalent integrable representations
of $(\PSU)^\circ$, and each irreducible
integrable \srep{} of \mbox{$(\PSU)^\circ$} is unitarily equivalent to
one on $M_N^\circ$.
Moreover, the integrable  *-representations
of $(\PSU)^\circ$ on $M_N^\circ$ and on
$P_N\hs \PS^{|N|+1}=  (((\PS^\circ)^{|N|+1})^\trans \circ (P_N)^\trans)^\trans$
are unitarily equivalent, where $\circ$ stands for the opposite multiplication.

With $n\in\N$, let $P\in\M_{n\times n}(\PS)$ be an invariant
projection. From \cite[Theorem 4.1]{SW2} and the preceding, we conclude that
the integrable representation of $(\PSU)^\circ$
on $P\hs \PS^n$ is equivalent to the direct sum of
irreducible representations on
$\big(\oplus_{N=-N_0}^{N_0} \oplus_{i=1}^{n_N(P)}\hs
P_N\big)\hs \PS^{m}$, where
$N_0\in\N$ is sufficiently large and
the orthogonal sum of projections is given as in Remark \ref{os}.
The numbers
$n_N(P)\in\N_0$ denote multiplicities and
$m=\sum_{N=-N_0}^{N_0} n_N(P) (|N|+1)$.
Set
$Q:=\oplus_{N=-N_0}^{N_0} \oplus_{i=1}^{n_N(P)}\hs P_N$
%$Q:=\oplus_{|N|<N_0}\,n_N(P)\,P_N$
and let $U:P\hs \PS^n\ra Q\hs \PS^m$ denote the unitary operator
realizing the equivalence.
Then the operator
$QUP\in \Hom_{\PS}(\C^n\otimes\PS ,\C^m\otimes\PS )$
 establishes a Murray-von Neumann equivalence between $P$ and $Q$.

Now let $P$ and $P'$ be invariant $\PS$-valued projections. The
proof of
\cite[Theorem 4.1]{SW2} shows that $P$ and $P'$ are
Murray-von Neumann equivalent if and only if $n_N(P)=n_N(P')$ for all
$N\in\Z$
in their decompositions described in the last paragraph.
From this, it follows first that $\KO^{\su}\hsp(\PS)$ has
cancellation, and second that $[P]-[P']=0$ if and only if
$n_N(P)=n_N(P')$ for all $N\in\Z$. Hence  $\KO^{\su}\hsp(\PS)$
is isomorphic to the free abelian group generated by $\{[P_N]:  N\in\Z\}$.
\end{proof}

Later, the index pairing in Proposition \ref{ipair} confirms that
(for transcendental $q$)
there are no relations between the generators $[P_N]$ of the
equivariant $\KO$-group $\KO^{\su}\hsp(\PS)$.

%*********************************************************
\section{Twisted Chern character }             \label{TCC}
%*********************************************************

The twisted Chern character will be defined as a map
from equivariant $\KO$-theory to twisted cyclic
homology. To do so, we need a convenient description of
twisted cyclic homology.

For a complex unital algebra $\A$ with automorphism $\lambda:\A\ra\A$,
set $C_n:=\A^{\otimes (n+1)}$ and define
\begin{align}                                                 \label{dni}
 d_{n,i} (a_0\otimes a_1\otimes \cdots \otimes a_n)
&= a_0\otimes \cdots \otimes a_i a_{i+1} \otimes
\cdots \otimes a_n, \quad i\neq n, \\
        d_{n,n}^\s (a_0\otimes a_1\otimes \cdots \otimes a_n)    \label{dnn}
&= \s(a_n) a_0\otimes a_1\otimes \cdots \otimes a_{n-1}, \\
\t_n^\s (a_0\otimes a_1\otimes \cdots \otimes a_n)
&=  \s(a_n) \otimes a_0 \otimes \cdots \otimes a_{n-1},         \label{tn}\\
        s_{n,i} (a_0\otimes a_1\otimes \cdots \otimes a_n)      \label{sni}
&= a_0\otimes\cdots a_i \otimes 1 \otimes a_{i+1} \otimes\cdots\otimes a_n,
\end{align}
where $n\in\N_0$ and $0\leq i\leq n$. For $ \lambda = \id$, these are
the face, cyclic and degeneracy operators of the
standard cyclic object associated to $\A$  \cite{L}.
For general $ \lambda $,  the operator $(\t_n^{\s})^{n+1}$
acting on $C_n$ by
$$
(\t_n^{\s})^{n+1} (a_0\otimes a_1\otimes \cdots \otimes a_n)
=  \s(a_0) \otimes \s(a_1) \otimes \cdots \otimes \s(a_{n})
$$
fails to be the identity.
To obtain a cyclic object, one  passes to the cokernels
$\Asnn := C_n / {\mathrm{im}}\hs (\id- (\t_n^\s)^{n+1})$.
The twisted cyclic homology
$\HC_\ast^\s (\A)$ is now the total homology of Connes' mixed
$(b,B)$-bicomplex $\BC^\s (\A)$:
\begin{equation}                                         \label{bB}
             \begin{CD}
        @ V{b_4^\s} VV @ V{b_3^\s} VV @ V{b_2^\s} VV @ V{b_1^\s} VV @ . @ . @ .\\
        {\Asnthree} @ <{B_2^\s}<< {\Asntwo} @ <{B_1^\s} << {\Asnone} @ <{B_0^\s} <<
        {\Asnzero} @ . @ . @ . @ .\\
         @ V{b_3^\s} VV @ V{b_2^\s} VV @ V{b_1^\s} VV @ . @ . @ . @ .\\
        {\Asntwo} @ <{B_1^\s} << {\Asnone} @ <{B_0^\s} << {\Asnzero} @ . @ . @ . @ . @ .\\
         @ V{b_2^\s} VV @ V{b_1^\s} VV @ . @ . @ . @ . @ .\\
        {\Asnone} @ <{B_0^\s} << {\Asnzero} @ . @ . @ . @ . @ . \\
         @ V{b_1^\s} VV @ . @ .  @ . @ . @ . @ .\\
        {\Asnzero} @ . @ . @ . @ . @ . @ .\\
            \end{CD}
\end{equation}
The boundary maps $b_n^\s$ and $B_n^\s$ are given by
\begin{equation}                                           \label{bn}
b_n^\s = \Sigma_{i=0}^{n-1} \; (-1)^i d_{n,i}+(-1)^n d_{n,n}^\s,\quad
        B_n^\s = (1- (-1)^{n+1} \t_{n+1} ) s_n N_n^\s,
\end{equation}
where
$N_n^\s = \Sigma_{j=0}^n \; (-1)^{nj} (\t_n^\s)^j$, and
the maps
\begin{equation}                                          \label{sn}
s_n  :  \Asnn \ra \Asnnplusone,\quad
s_n ( a_0\otimes a_1\otimes \cdots \otimes a_n) =
1\otimes a_0\otimes a_1\otimes \cdots \otimes a_n
\end{equation}
are called ``extra degeneracies''.

The homology of the columns is the twisted Hochschild homology of $\A$
and denoted by $\HH^\s_\ast(\A)$.
An element $\eta\in \A^{\otimes (n+1)}$ such that $b_n^\s (\eta)=0$ is called
a {\it twisted Hochschild $n$-cycle}.
The class of $\eta$ in $C^\s_n$ defines an element in
twisted Hochschild  and cyclic homology by putting it at the $(0,n)$-th position
of the $(b,B)$-bicomplex \eqref{bB}.

By dualizing, one passes from twisted cyclic homology $\HC_\ast^\s (\A)$
to twisted cyclic cohomology $\HC^\ast_\s (\A)$. A complex for
computing $\HC^\ast_\s (\A)$ is obtained by applying the functor
$\Hom(\hs \cdot\hs,\C)$ to each entry of the $(b,B)$-complex.

An alternative description of twisted cyclic cohomology $\HC^\ast_\s (\A)$
is as follows \cite{KMT}.
Let $\ASN$ denote the space of of $(n+1)$-linear forms $\phi$ on $\A$
such that $\phi=(-1)^n \t^{\s\ast}_n \phi$, where
$$
\t^{\s\ast}_n \phi(a_0,\mdots,a_n)=\phi(\s(a_n),a_0,\mdots,a_{n-1}).
$$
With the coboundary operator $b^{\s\ast}_n : \ASNN\ra\ASN$,
\begin{align*}
b^{\s\ast}_n\phi(a_0,\mdots,a_n)
=\sum_{j=0}^{n-1}(-1)^j\phi(a_0,\mdots, a_j & a_{j+1},\mdots,a_n)\\
&+(-1)^n\phi(\s(a_n)a_0,\mdots,a_{n-1}),
\end{align*}
one gets a cochain complex whose homology is isomorphic to
$\HC^\ast_\s (\A)$.
The elements $\phi\in \ASN$ satisfying $b^{\s\ast}_{n+1}\phi=0$
are called \emph{twisted cyclic $n$-cocycles}.
An isomorphism between the two versions of twisted cyclic cohomology
is given by putting a twisted cyclic $n$-cocycle at the
$(0,n)$-th position in the dual $(b,B)$-complex and zeros elsewhere.

Evaluating cycles on cocycles yields a dual pairing between
$\HC_\ast^\s (\A)$ and $\HC^\ast_\s (\A)$.
For $\s=\id$, there is a Chern character map from $K$-theory to
cyclic homology available. Composing the Chern character
with the evaluation on cocycles defines a pairing between
$K$-theory and cyclic cohomology (the Chern-Connes pairing).
Our aim is to construct a similar  pairing between
equivariant $\KO$-theory and twisted cyclic cohomology.
The primary tool will be a ``twisted'' Chern character
from the equivariant $\KO$-group to even twisted cyclic homology.
This shall be our concern in the remainder of this section.

The discussion will be restricted to the following setting:
We assume that $\U$ is a Hopf *-algebra, $\B$ a unital left
$\U$-module *-algebra, and the automorphism $\s:\B\ra\B$ can be described
by a group-like element $k\in\U$, i.e., $\Delta(k)=k\otimes k$
and $\s(b)=k\rt b$ for all $b\in \B$. Note that $\vare(k)=1$ and $S(k)=k^{-1}$.
In particular, it follows that  $\s(b)^*=\s^{-1}(b^*)$.

The link between the ``non-twisted'' and the ``twisted'' case is provided
by the so-called quantum trace. Given  matrices $A_k=(a^k_{i_k,j_k})\in\M_{m\times m}(\B)$,
and  an \mbox{(anti-)}representation $\rho^\circ:\U\ra\End(\C^m)$,
we define the quantum trace $\tr_\s$ by
\begin{equation}                                              \label{trq}
 \tr_\s(A_0\otimes  A_1\otimes \cdots\otimes A_n)= \msum_{j_0,\ldots,
j_{n+1}} \;
\rho^\circ(k)_{j_{n+1},j_0}\hs
a^0_{j_0,j_1}\otimes a^1_{j_1,j_2} \otimes \cdots \otimes a^n_{j_n,j_{n+1}}.
\end{equation}
\begin{lem}                                                  \label{Trq}
With the conventions introduced above, let $\A$ denote the
algebra of  $\adro$-invariant matrices belonging to  $\M_{m\times m}(\B)$.
Then $\tr_\s$ defines a morphism of complexes from
$\BC^\id (\A)$ to $\BC^\s (\B)$.
\end{lem}
\begin{proof}
To prove the lemma, it suffices to show that $\tr_\s$ intertwines the operators
defined in Equations \eqref{dni}--\eqref{sni} and \eqref{sn}. For the face
operators  $d_{n,i}$, $i\neq n$, and the degeneracy operators $s_{n,i}$ and
$s_n$, the assertion is obvious. For $d_{n,n}^\s$ and $\t_n^\s$, one uses the fact
that
$$
    \rho^\circ(k)\hs A =\s(A)\hs \rho^\circ(k)\quad        \mbox{for all}\ A\in\A,
$$
since
$$
A=\adro(k)(A)= \rho^\circ(k^{-1}) \hs(k\rt A)\hs \rho^\circ(k)=
\rho^\circ(k^{-1})\hs \s(A)\hs\rho^\circ(k)
$$
for any  $\adro$-invariant matrix $A\in\M_{m\times m}(\B)$.
\end{proof}

The next proposition introduces a twisted version of the Chern character
mapping equivariant $\KO$-groups into twisted cyclic homology.

\begin{prop}                                            \label{chern}
Let $\U$, $\B$, $\s$, $\rho^\circ$ be as above. Suppose that there is an
automorphism of $\B$ satisfying the conditions of Definition \ref{MvN}.
For any invariant projection $P\in\M_{m\times m}(\B)$, set
\begin{equation*}
\ch_{2n}(P):=(-1)^n\, \mbox{$\frac{(2n)!}{n!}$} \,\tr_\s
(P\otimes P\otimes \cdots \otimes P)\,\in \,\B^{\otimes 2n+1}.
\end{equation*}
Then there are well-defined additive maps
$
\chn :\KO^{\U}\hsp(\B)\ra \HC^\s_{2n}(\B)
$
given by
$$
\chn([P])=(\ch_{2n}(P),\mdots, \ch_{0}(P)).
$$
\end{prop}
\begin{proof}
Let $P\in \M_{m\times m}(\B)$ be an invariant projection.
From the Chern character of $\KO$-theory with values in (non-twisted) cyclic
homology (cf.\ \cite{L}), it is known that
$(\chid_{2n}(P),\mdots, \chid_{0}(P))$ defines a cycle
in $\BC^\id (\A)$,
where $\A$ denotes again the subalgebra of  $\adro$-invariant matrices
in $\M_{m\times m}(\B)$. By Lemma \ref{Trq},
$(\ch_{2n}(P),\mdots, \ch_{0}(P))$ is a  cycle in
$\BC^\s (\B)$.
Moreover, the additivity of
$\chn$ follows from the additivity of $\tr_\s$.

It remains to prove that the homology class of $\chn([P])$
does not depend on representatives.
Let $P\in \M_{j\times j}(\B)$ and $Q\in \M_{k\times k}(\B)$ be
invariant projections with respect to the anti-representations
$\rho_1^\circ:\U\ra\End(\C^j)$ and $\rho_2^\circ:\U\ra\End(\C^k)$,
respectively,
and suppose that $P$ and $Q$ are Murray-von Neumann equivalent.
By considering the direct sum
$\rho_1^\circ\oplus\rho_2^\circ$ on $\C^{j+k}$, we may assume that
$P$ and $Q$ belong to the
same matrix algebra $\M_{m\times m}(\B)$ and that there is an invertible
element $U\in\M_{m\times m}(\B)$ such that $UPU^{-1}=Q$.
To be more precise, we consider
$$
\left( \begin{matrix} P & 0\\
0 & 0 \end{matrix}\right)\sim P, \quad
\left( \begin{matrix} 0 & 0\\
0 & Q \end{matrix}\right)\sim Q, \quad
U:=\left( \begin{matrix} 1-P  & PV^\dagger Q\\
QVP&1-Q \end{matrix}\right)=U^{-1},
$$
where $V\in \Hom_\B(\C^j\otimes\B ,\C^k\otimes\B )$
establishes the Murray-von Neumann equivalence between $P$ and $Q$.
Clearly, $P$ and $Q$ are $\adro$-invariant with respect to
$\rho_1^\circ\oplus\rho_2^\circ$ and, by Definition \ref{MvN},
so is $U$. Since conjugation acts as the identity on cyclic homology
\cite[Proposition 4.1.2]{L}, $(\chid_{2n}(P),\mdots, \chid_{0}(P))$ and
$(\chid_{2n}(Q),\mdots, \chid_{0}(Q))$ differ in $\BC^\id (\A)$ only by
a boundary. Hence,
by Lemma \ref{Trq}, the Chern characters
$\chn([P])$ in $\HC^\s_{2n}(\B)$ do not depend on the
representative of the class $[P]$ in  $\KO^{\U}\hsp(\B)$.
\end{proof}
\begin{cor}                                           \label{pair}
There is a pairing
$$
\ip{\cdot}{\cdot}:\HC_\s^{2n}(\B)\times \KO^{\U}\hsp(\B)\ra\C
$$
given by evaluating cocycles on $\chn$.
For a twisted cyclic $2n$-cocycle $\Phi$,
it takes the form
\begin{equation}                                               \label{chpair}
\ip{[\Phi]}{[P]}=(-1)^n\, \mbox{$\frac{(2n)!}{n!}$} \,
\Phi\big(\tr_\s(P\otimes P\otimes \cdots \otimes P)\big).
\end{equation}
\end{cor}
\begin{proof}
The first assertion follows from Proposition \ref{chern}, the second from the
embedding of cycles $\Phi\in \ASN$ satisfying $b^{\s\ast}_{n+1}\Phi=0$ into
the $(b,B)$-complex.
\end{proof}

\begin{rem}                                                         \label{RemS}
Note that the pairing between $\Phi$ and Chern character $\ch_{2n}$
is compatible with Connes' periodicity operator
$\cS:\HC_\s^{n}(\B)\ra \HC_\s^{n+2}(\B)$, that is,
\begin{equation}
\cS\Phi(\ch_{2n+2}(P))=\Phi(\ch_{2n}(P))                           \label{Sop}
\end{equation}
for any invariant idempotent $P$, where
\begin{align*}
\cS\Phi(a_0,\ldots,a_{n+2}):=&-\mbox{$\frac{1}{(n+1)(n+2)}$}\Big(
\sum_{i=1}^{n+1}\phi(a_0,\ldots,a_{i-1}a_ia_{i+1},\ldots,a_{n+2})\\
&-\!\sum_{1\leq i<j\leq (n+1)}\hspace{-17pt}(-1)^{i+j}
\phi(a_0,\ldots,a_{i-1}a_i,\ldots, a_ja_{j+1},\ldots,a_{n+2})\Big).
\end{align*}
Actually, Corollary \ref{pair} defines a pairing between periodic twisted
cyclic cohomology and equivariant $\KO$-theory, but we shall not go
into the details.
\end{rem}

\begin{rem}                                                       \label{tp}
The description of $\HC^\ast_\s (\B)$ involves only an algebra $\B$ with an
automorphism $\s$. If $\B$ is a *-algebra and $\s$ satisfies
$\s(a)^*=\s^{-1}(a^*)$, then, by Example \ref{Ks},  we have a definition of
equivariant $\KO$-theory  $\KO^\s\hsp(\B)$ and, by Corollary \ref{pair}, a
pairing between $\HC^\ast_\s (\B)$ and  $\KO^\s\hsp(\B)$.
\end{rem}

%*********************************************************
\section{Chern-Connes pairing}                         \label{CCP}
%*********************************************************

%+++++++++++++++++++++++++++++++++++++++++++++++++++++++++++
\subsection{The modular automorphism case}              \label{mac}
%+++++++++++++++++++++++++++++++++++++++++++++++++++++++++++

This section is devoted to the calculation of the pairing
\begin{equation}                                               \label{dispair}
\ip{\cdot}{\cdot}:\HC_\theta^{2n}(\PS)\times \KO^{\su}\hsp(\PS)\ra\C
\end{equation}
described in Corollary \ref{pair},
where $\theta$ denotes the modular automorphism from Equation \eqref{Phaut}.
This includes in particular the pairing of $\KO^{\su}\hsp(\PS)$ with the
$\theta$-twisted cyclic 2-cocycle $\tau$
associated with the volume form of
the distinguished 2-dimensional covariant differential calculus on $\PS$
(cf.\ Proposition \ref{cycle}).
Note, by the way, that
$\KO^{\su}\hsp(\PS)$ is given in Proposition \ref{eqK0}
with respect to the inverse automorphism
$\theta^{-1}$.

The twisted cyclic homology of $\PS$ was computed by Hadfield \cite{H}. Dualizing, we
conclude from \cite{H} that
$$
\HC_\theta^{2n}(\PS)=\C[\cS^n\vare]\oplus \C[\cS^n h],\quad
\HC_\theta^{2n+1}(\PS)=0, \quad n\in\N_0,
$$
where $h$ is the Haar state on $\PS$, $\vare$ is the restriction
of the counit of $\SU$ to $\PS$, and $\cS:\HC_\theta^{n}(\PS)\ra \HC_\theta^{n+2}(\PS)$
denotes Connes' periodicity operator from Remark \ref{RemS}.
As a consequence, the pairing in \eqref{dispair} is completely determined
by the pairing of $\HC_\theta^{0}(\PS)$ with $\KO^{\su}\hsp(\PS)$.

\begin{prop}                                              \label{ipair}
The pairing between $\HC_\theta^{0}(\PS)$ and $\KO^{\su}\hsp(\PS)$
is given by
\begin{align*}
\ip{[\vare]}{[P_N]}= q^{N},\quad \ip{[h]}{[P_N]}= q^{-N},\quad N\in\Z.
\end{align*}
\end{prop}

\begin{proof}
Set $n:=N/2$ and $l:=|N|/2$.
With the notation of Section \ref{sec-equiv},
it follows from \eqref{ft} and \eqref{rhoo} that
$\fett{t}^{l}_n\hs\rho_l^\circ(K^{-2})=K^{2}\rt \fett{t}^{l}_n$.
Applying successively
\eqref{chpair}, \eqref{PN}, \eqref{trq}, \eqref{EFK} and \eqref{hopf},
we get
\begin{align*}
\ip{[\vare]}{[P_N]}&=\vare(\cht_0(P_N))=
\vare\big( \tr \,\rho_{l}^\circ(K^{-3})\,
\fett{t}^{l*}_{n}\,\fett{t}^{l}_{n}\,\rho_{l}^\circ(K)\big)
=\vare\big( \tr \,
\fett{t}^{l*}_{n}\,(K^2\rt \fett{t}^{l}_{n})\big)\\
&=\mbox{$\sum_{j=-l}^l$}   q^{2j} \vare( t^{l*}_{n,j}t^{l}_{n,j})
=q^{2n}.
\end{align*}
Note that $\theta(t^{l}_{n,j})=q^{-2j-2n}t^{l}_{n,j}$ by \eqref{haut} and
$\sum_{j=-l}^l t^{l}_{n,j} t^{l*}_{n,j}=1$ by \eqref{hopf}.
Analogously to the above, we have
\begin{equation*}
\ip{[h]}{[P_N]}=\mbox{$\sum_{j=-l}^l$}   q^{2j} h (t^{l*}_{n,j}t^{l}_{n,j})
= q^{-2n}  \mbox{$\sum_{j=-l}^l$}   h (t^{l}_{n,j} t^{l*}_{n,j}) =q^{-2n}.
\\[-24pt]
\end{equation*}
\end{proof}\smallskip

Now we compute the pairing of $\KO^{\su}\hsp(\PS)$ with the $\theta$-twisted cyclic
2-cocycle $\tau$ from Proposition \ref{cycle}.

\begin{prop}
The pairing of $[\tau]$ with
$\KO^{\su}\hsp(\PS)$  is given by
\[                                                     \label{iptau}
\ip{[\tau]}{[P_N]}= 2[N]_q,\quad N\in\Z.
\]
In particular, $\tau=\frac{2}{q^{-1}-q} (\cS h - \cS \vare)$.
\end{prop}

\begin{proof}
Let $n=N/2>0$.
Inserting \eqref{cyc} and \eqref{PN} into \eqref{chpair}, we get
$$
\ip{[\tau]}{[P_N]}=2\hs h(\tr_\theta \,
\fett{t}^{n*}_{n}\,\fett{t}^{n}_{n}
\big(q^{-1}\dif_F(\fett{t}^{n*}_{n}\,\fett{t}^{n}_{n})
\dif_E(\fett{t}^{n*}_{n}\,\fett{t}^{n}_{n})-
q\dif_E(\fett{t}^{n*}_{n}\,\fett{t}^{n}_{n})
\dif_F(\fett{t}^{n*}_{n}\,\fett{t}^{n}_{n})
\big)).
$$
From \eqref{v}, \eqref{dif} and \eqref{difab}, it follows that
\begin{equation}                                                    \label{dift}
\dif_E(\fett{t}^{n*}_{n}\,\fett{t}^{n}_{n})
=-q^{n-1}\alpha^n_{n-1}\fett{t}^{n*}_{n}\,\fett{t}^{n}_{n-1},\quad
\dif_F(\fett{t}^{n*}_{n}\,\fett{t}^{n}_{n})
=q^{n}\alpha^n_{n-1}\fett{t}^{n*}_{n-1}\,\fett{t}^{n}_{n}.
\end{equation}
Equation \eqref{t*t} implies
$\fett{t}^{n}_{n}\fett{t}^{n*}_{n}=\fett{t}^{n}_{n-1}\fett{t}^{n*}_{n-1}=1$
and $\fett{t}^{n}_{n}\fett{t}^{n*}_{n-1}=0$. Hence
$$
\ip{[\tau]}{[P_N]}=2q^{2n}(\alpha^n_{n-1})^2
\hs h(\tr_\theta \,\fett{t}^{n*}_{n}\,\fett{t}^{n}_{n})=2q^{2n}(\alpha^n_{n-1})^2
\hs h(\tr_\theta \,P_N).
$$
Applying $h(\tr_\theta \,P_N)=\ip{[h]}{[P_N]}=q^{-2n}$
and $(\alpha^n_{n-1})^2=[2n]_q$ proves \eqref{iptau} for $N>0$.
If $n=N/2<0$, Equation \eqref{dift} becomes
\[                                                  \tag{\ref{dift}'}
\dif_E(\fett{t}^{|n|*}_{n}\,\fett{t}^{|n|}_{n})
=q^{n}\alpha^{|n|}_{n}\fett{t}^{|n|*}_{n+1}\,\fett{t}^{|n|}_{n},\quad
\dif_F(\fett{t}^{|n|*}_{n}\,\fett{t}^{|n|}_{n})
=-q^{n+1}\alpha^{|n|}_{n}\fett{t}^{|n|*}_{n}\,\fett{t}^{|n|}_{n+1},
\]
and by the same arguments as above, we get
$$
\ip{[\tau]}{[P_N]}=-2q^{2n}(\alpha^{|n|}_{n})^2
\hs  \ip{[h]}{[P_N]} =-2\hs[2|n|]_q=2[N]_q.
$$
The case $N=0$ is trivial since $\dif_F(1)=\dif_E(1)=0$.

It has been shown in \cite{H} that
$\tau=\beta (\cS h - \cS \vare)$ with $\beta\in\C\setminus \{0\}$.
By the preceding, Equation \eqref{Sop} and  Proposition \ref{ipair},
\begin{align*}
2[N]_q&=\tau(\cht_2(P_N))=\beta\hs \cS( h -  \vare)(\cht_2(P^N_N))
=\beta\hs  ( h - \vare)(\cht_0(P^N_N))\\
&=\beta \hs (q^{-N}-q^N),
\end{align*}
thus $\beta =\frac{2}{q^{-1}-q}$.
\end{proof}

%+++++++++++++++++++++++++++++++++++++++++++++++++++++++++++
\subsection{The case of the inverse modular automorphism}   \label{imac}
%+++++++++++++++++++++++++++++++++++++++++++++++++++++++++++

The modular automorphism $\theta$ does not correspond to the
``no dimension drop'' case of twisted Hochschild homology since
$\HH_2^{\theta}(\PS)=0$ (see \cite{H}). On the other hand,
as shown in \cite{H}, the dimension drop can be avoided by
taking the inverse modular automorphism. In this case,
$\HH_2^{\theta^{-1}}(\PS)\cong\C$ and
$\HC_2^{\theta^{-1}}(\PS)\cong\C^2$.
By duality, we have $\HC^2_{\theta^{-1}}(\PS)\cong\C^2$.

Observe that any $\s$-cyclic $n$-cocycle $\psi\in \HC^n_{\s}(\A)$
defines a functional on $\HH_n^{\s}(\A)$ since
$\psi\circ b_{n+1}^\s= b_{n+1}^{\s *}\psi =0$ (cf.\ Section \ref{TCC}).
One generator of $\HC^2_{\theta^{-1}}(\PS)$ is given by
$[\cS\vare]$. It is easily shown by replacing $\tr_\theta$ by
$\tr_{\theta^{-1}}$ in the previous section that
$\ip{[\cS\vare]}{[P_N]}=\vare(\mathrm{ch}^{\theta^{-1}}_0(P_N))= q^{-N}$
for all $N\in\Z$. This twisted 2-cocycle descends to the trivial
functional on $\HH_2^{\theta^{-1}}(\PS)$.

The other generator of $\HC^2_{\theta^{-1}}(\PS)$,
which we denote by $[\phi]$,
was recently described
by Kr\"ahmer \cite{Kr} and truly corresponds to the
``no dimension drop'' case in the sense that it is non-trivial on
$\HH_2^{\theta^{-1}}(\PS)$. In the following, we
compute the pairing of $\KO^{\su}\hsp(\PS)$ with $[\phi]$.

It has been shown in \cite{Kr} that
$[\phi]=[\varphi + b^{\theta^{-1}*}_2 \chi]$,
where
\begin{equation}                                   \label{coc}
\varphi(a_0,a_2,a_3)=q\hs\vare(a_0 (K^{-1}E\rt a_1)(K^{-1}F\rt a_1)),\quad
a_0,a_1,a_2\in\PS,
\end{equation}
and $\chi: \PS\otimes\PS\ra\C$ is a linear functional
such that
\begin{equation}                                   \label{phi}
(\varphi + b^{\theta^{-1}*}_2 \chi)(1,a_1,a_2)=0\quad
\mbox{for all}\ \,a_1,a_2\in\PS.
\end{equation}
Using the vector space basis
$\{A^lB^k,\ A^mB^{n*}: l,k,m\in\N_0,\ n\in\N\}$ of $\PS$,
$\chi$ is determined by
\begin{equation}                                  \label{chi}
\chi(1,A)=(q^{-1}-q)^{-1},\quad \chi(1,A^2)=(q-q^3)^{-1}, \quad
\chi(A,A)=(2(q-q^3))^{-1}
\end{equation}
and $\chi(A^lB^k,\ A^mB^{n*})=0$ otherwise.

Observe that the evaluation of $\phi$ on Hochschild cycles
reduces to the application of $\varphi$ from Equation \eqref{coc}.
Therefore we will slightly change the Chern character mapping
$\mathrm{ch}^{\theta^{-1}}_2$ in Proposition \ref{chern}
in order to obtain Hochschild 2-cycles.
\begin{prop}                                                     \label{ochi}
For the projections $P_N$ defined in Equation \eqref{PN}, set
\begin{equation}                                                  \label{och}
\och(P_N):=\tr_\thetai\big((1-2P_N)\otimes(P_N-\vare(P_N))\otimes(P_N-\vare(P_N))\big),
\end{equation}
where the counit $\vare$ of $\SU$ is applied to each
component of $P_{N}$.
Then $\och(P_N)$ is a $\thetai$-twisted Hochschild 2-cycle, i.e.,
$b_2^\thetai( \och(P_N))=0$.
\end{prop}
\begin{proof}
Note that $\vare(P_N)$ is a diagonal complex matrix. In particular, it
commutes with $\rho_{|N|/2}^\circ(K^{2})$.
Using $\tr_\thetai\vare(P_N)X\otimes Y=\tr_\thetai X\otimes Y\vare(P_N)$ and
$\tr_\thetai X\vare(P_N)\otimes Y=\tr_\thetai X\otimes\vare(P_N) Y$ for any
$X,Y\in \M_{|N|+1\times |N|+1}(\PS)$, one easily shows that
\begin{equation}
b_2^\thetai( \och(P_N))=\tr_\thetai\big( 1\otimes (P_N-\vare(P_N))\big)
=1\otimes (\tr_\thetai P_N -\tr_\thetai \vare(P_N)).
\end{equation}
Thus it remains to show that $\tr_\thetai P_N=\tr_\thetai \vare(P_N)$. Let $l:=N/2$.
Then
\begin{align*}
F\rt \tr_\thetai P_N&= \sum{}_{k=-|l|}^{|l|}\,q^{-2k}\big(
(F\rt t^{|l|*}_{lk})(K\rt t^{|l|}_{lk})+(K^{-1}\rt t^{|l|*}_{lk})(F\rt t^{|l|}_{lk})\big)\\
&=\sum{}_{k=-|l|}^{|l|}\, -q^{-(k+1) }\alpha_k^{l} t^{|l|*}_{l,k+1}t^{|l|}_{lk}
+ q^{-k}\alpha^{l}_{k-1} t^{|l|*}_{l,k}t^{|l|}_{l,k-1} \,=\,0
\end{align*}
by Equations \eqref{lmod} and \eqref{EFK}.
Similarly, $K\rt \tr_\thetai P_N=\tr_\thetai P_N$, hence $\tr_\thetai P_N$
belongs to the spin 0 representation of  $M_0=\PS$. Therefore there exists an
$\alpha\in\C$ such that $\tr_\thetai P_N=\alpha\hs t^0_{00}$.  Since  $t^0_{00}=1$,
we can write (with a slight abuse of notation)
$\tr_\thetai \vare(P_N)=\vare(\tr_\thetai P_N)=\alpha=\tr_\thetai P_N$
which concludes the proof.
\end{proof}
\begin{rem}
Actually, $\och$ and its $2n$-dimensional analogs give rise to a
twisted Chern character
from equivariant $\KO$-theory into twisted cyclic homology if we replace the
$(b,B)$-bicomplex \eqref{bB} by the so-called
reduced $(b,B)$-bicomplex. The reduced $(b,B)$-bicomplex is defined
analogous to the $(b,B)$-bicomplex \eqref{bB} but with
$C_n:=\A\otimes \bar\A^{\otimes (n)}$, where $\bar\A:=\A/\C$, see
\cite{L} and \cite{GFV} for details in the ``non-twisted'' case.
\end{rem}

We turn now to the computation of the Chern-Connes pairing in the
``no dimension drop'' case.
\begin{prop}
The pairing of $[\phi]$ with
$\KO^{\su}\hsp(\PS)$ yields
$$
\ip{[\phi]}{[P_N]}=  [N]_q,\quad N\in\Z.
$$
\end{prop}
\begin{proof}
We first claim that
$\phi(\mathrm{ch}^{\theta^{-1}}_2\!(P_N))=\phi(\och(P_N))$.
By Equation \eqref{phi},
$\phi(1, a, b)=0$ for all $a,b\in\PS$.
Since $\vare(P_N)\in \M_{|N|+1\times |N|+1}(\C)$, it suffices to show that
$\phi(a,1, b)=\phi(a, b, 1)=0$ for all $a,b\in\PS$.
For $\varphi$, this follows from $E\rt 1=F\rt 1=0$.
Observe that, by \eqref{chi}, $\chi(x,1)=0$ for all $x\in\PS$.
Hence $b^{\theta^{-1}*}_2 \chi(a,1,b)=\chi(\theta^{-1}(b)a,1)=0$ and
$b^{\theta^{-1}*}_2 \chi(a,b,1)=\chi(ab,1)=0$
which proves the claim.

Since $b_2^\thetai( \och(P_N))=0$, it follows that
$\ip{[\phi]}{[P_N]}= \varphi(\och\hsp(P_N) )$.
Set $l:=|N|/2$ and $n:=N/2$. Inserting \eqref{PN} and \eqref{och}
into \eqref{coc} gives
$$
\ip{[\phi]}{[P_N]}=q\hs\tr\, \sigma_l(K^{-2} )\hs
\vare(1-2 \fett{t}^{l*}_{n}\,\fett{t}^{l}_{n})\hs
\vare(K^{-1}E\rt \fett{t}^{l*}_{n}\,\fett{t}^{l}_{n})\hs
\vare(K^{-1}F\rt \fett{t}^{l*}_{n}\,\fett{t}^{l}_{n}).
$$
By \eqref{lmod} and \eqref{ft},
$K^{-1}E\rt \fett{t}^{l*}_{n}\,\fett{t}^{l}_{n}
=\sigma_l(-qEK)\fett{t}^{l*}_{n}\,\fett{t}^{l}_{n}
+\sigma_l(K^2)\fett{t}^{l*}_{n}\,\fett{t}^{l}_{n}\sigma_l(K^{-1}E)$
and
$K^{-1}F\rt \fett{t}^{l*}_{n}\,\fett{t}^{l}_{n}
=\sigma_l(-q^{-1}FK)\fett{t}^{l*}_{n}\,\fett{t}^{l}_{n}
+\sigma_l(K^2)\fett{t}^{l*}_{n}\,\fett{t}^{l}_{n}\sigma_l(K^{-1}F)$.
The complex matrices $\sigma_l(f)$, $f\in\su$, can be derived from
the formulas in \eqref{EFK}.
Inserting
$\vare(\fett{t}^{l*}_{n}\,\fett{t}^{l}_{n})=(\delta_{ni}\delta_{nj})_{i,j=-l}^l$
yields
$\sigma_l(K^{-2} )\hs \vare(1-2 \fett{t}^{l*}_{n}\,\fett{t}^{l}_{n})
=\big((-1)^{\delta_{nj}}q^{-2j}\delta_{ij}\big)_{i,j=-l}^l$.
If  $N\geq 0$, then $n=l$, thus
$\vare(K^{-1}E\rt \fett{t}^{l*}_{n}\,\fett{t}^{l}_{n})
=\big(q^{l}[2l]_q^{1/2}\delta_{li} \delta_{l-1,j}  \big)_{i,j=-l}^l$
and
$\vare(K^{-1}F\rt \fett{t}^{l*}_{n}\,\fett{t}^{l}_{n})
=\big(-q^{l-1}[2l]_q^{1/2}\delta_{l-1,i} \delta_{lj}  \big)_{i,j=-l}^l$
by the above.
Multiplying the matrices and taking the trace gives
$\ip{[\phi]}{[P_N]}=[2l]_q=[N]_q$. A similar calculation shows that
$\ip{[\phi]}{[P_N]}=-[2l]_q=[N]_q$ for $N<0$.
\end{proof}
%
%
%*********************************************************
\section{Orientation}                           \label{O}
%*********************************************************

In this section, we show that there exists a twisted Hochschild
2-cycle $\eta$ in $\PS^{\otimes 3}$ such that $\pi_D(\eta)=\gamma_q$, where
$$
\pi_D(\mbox{$\sum_j$} a^0_j\otimes a^1_j\otimes a^2_j):=
\mbox{$\sum_j$} a^0_j[D,a^1_j][D,a^2_j]\quad \mbox{and} \quad
\gamma_q:= \left(
\begin{matrix} q^{-1} & 0 \\ 0 & -q \end{matrix} \right).
$$
In analogy to axiom $(4')$ in \cite{C2}, we call $\eta$ a choice of orientation.

In the commutative case,
the Hochschild cycle corresponds to the volume form which is non-trivial
in de Rham cohomology. For this reason, we consider  again the inverse
modular automorphism $\theta^{-1}$ which avoids the dimension drop.
\begin{prop}                                                \label{ori}
With $P_1$ given in Equation \eqref{PN} for $N=1$, define
$$
\eta:=\och(P_1)=\tr_{\theta^{-1}}\big(
(1-2P_{1})\otimes ( P_{1}-\vare(P_{1}))
\otimes ( P_{1}-\vare(P_{1}))\big).
$$
Then $b_2^{\theta^{-1}}(\eta)=0$,
the class of $\eta$ in $\HH_2^{\theta^{-1}}(\PS)$ and
$\HC_2^{\theta^{-1}}(\PS)$ is non-zero,
$$
\HH_2^{\theta^{-1}}(\PS)\cong \C [\eta],\quad
\HC_2^{\theta^{-1}}(\PS)\cong \C [1]\oplus \C [\eta],
$$
and $\pi_D(\eta)=\gamma_q$.
\end{prop}
\begin{proof}
Since $\eta= \och(P_{1})$, it follows immediately from Proposition \ref{ochi}
that $\eta$ is a $\theta^{-1}$-twisted Hochschild 2-cycle.
For brevity of notation, set $n:=1/2$. By \eqref{Da} and \eqref{PN},
$$
\pi_D(\eta)=
\left(\begin{matrix}\tr_{\theta^{-1}}\hs(1-2 \fett{t}^{n*}_{n}\,\fett{t}^{n}_{n})
\dif_E(\fett{t}^{n*}_{n}\,\fett{t}^{n}_{n})
\dif_F(\fett{t}^{n*}_{n}\,\fett{t}^{n}_{n}) &\mbox{ }\hspace{-20pt} 0  \\
\mbox{ }\hspace{-60pt} 0 & \hspace{-80pt}
\tr_{\theta^{-1}}\hs(1-2 \fett{t}^{n*}_{n}\,\fett{t}^{n}_{n})
\dif_F(\fett{t}^{n*}_{n}\,\fett{t}^{n}_{n})
\dif_E(\fett{t}^{n*}_{n}\,\fett{t}^{n}_{n})
 \end{matrix} \right).
$$
From Equation \eqref{dift} and the argument following it, we conclude that
$$
\pi_D(\eta)=
\left(\begin{matrix}\tr\,\rho_n^\circ(K^2)\,\fett{t}^{n*}_{n}\,\fett{t}^{n}_{n} & 0\\
0 & - \tr\,\rho_n^\circ(K^2)\,\fett{t}^{n*}_{-n}\,\fett{t}^{n}_{-n}
\end{matrix} \right).
$$
Straightforward computations using the embedding \eqref{ABB*} and the relations
in $\SU$ (cf.\ \cite{KS}) show that
$$
\fett{t}^{n*}_{n}\,\fett{t}^{n}_{n}
=\left(\begin{matrix} A & B^*\\ B & 1-q^2A\end{matrix} \right),~~
\fett{t}^{n*}_{-n}\,\fett{t}^{n}_{-n}
=\left(\begin{matrix}  1-A & -B^*\\ -B & q^2A\end{matrix} \right), ~~
\rho_n^\circ(K^2)=\left(\begin{matrix} q & 0\\ 0 & q^{-1} \end{matrix} \right).
$$
Inserting these matrices into the previous equation and taking the trace
gives $\pi_D(\eta)=\gamma_q$.

Observe that
\begin{align}                                                        \label{etat}
\eta &= \tr\,\rho_n^\circ(K^2) \hs
(1\hsp-\hsp 2 \fett{t}^{n*}_{n}\,\fett{t}^{n}_{n})\otimes
( \fett{t}^{n*}_{n}\,\fett{t}^{n}_{n}\hsp-\hsp\vare(\fett{t}^{n*}_{n}\,\fett{t}^{n}_{n}))
\otimes(\fett{t}^{n*}_{n}\,\fett{t}^{n}_{n}\hsp-\hsp\vare(\fett{t}^{n*}_{n}\,\fett{t}^{n}_{n}))\\
\notag &= q(1\hsp -\hsp 2A)\otimes(A\otimes A+ \B^*\otimes B)
-2q B^*\otimes (B\otimes A-q^2A\otimes B) \\
\notag & -2q^{-1}B\otimes (A\otimes B^*\hsp  -\hsp q^2 B^*\otimes A)
+q^{-1}(2q^2A\hsp -\hsp 1)\otimes(B\otimes B^*\hsp  +\hsp  q^4 A\otimes A).
\end{align}
Let $\omega_2$ denote the $\theta^{-1}$-twisted 2-cycle
defined in \cite[Equation (27)]{H}. Then
$$
\omega_2-q^{-1}\eta = 2(B\otimes B^*\otimes A - q^{-2}B\otimes A \otimes B^*
 - q^2  B^*\otimes A\otimes B + B^*\otimes B\otimes A).
$$
Setting
\begin{align*}
\beta&=q^{-2} 1\otimes B\otimes B^*\otimes A - q^2 1\otimes  B^*\otimes B\otimes A
+ q^{-2} 1\otimes A \otimes B\otimes B^*\\
&\qquad \quad - q^2 1\otimes  A \otimes B^*\otimes B
+ q^4 1\otimes  B^*\otimes A\otimes B -q^{-4}  1\otimes B\otimes A \otimes B^*
\end{align*}
one easily checks that
$b_3^{\theta^{-1}}(2(q^{-2}-q^2)^{-1}\beta)= \omega_2-q^{-1}\eta$,
hence $[\eta]=[q\hs\omega_2]$ in $\HH_2^{\theta^{-1}}\hsp (\PS)$
and  $\HC_2^{\theta^{-1}}\hsp (\PS)$. The last statements of Proposition \ref{ori}
follow now from the results in \cite{H}.
\end{proof}

The next proposition shows that $\eta$ represents the
volume form of the covariant differential calculus
on $\PS$.
\begin{prop}
For $\mbox{$\sum_j$} a^0_j\otimes a^1_j\otimes a^2_j\in \PS^{\otimes 3}$, let
$$
\pi_\wedge(\mbox{$\sum_j$} a^0_j\otimes a^1_j\otimes a^2_j):=
\mbox{$\sum_j$} a^0_j\, \dd a^1_j\wedge \dd a^2_j.
$$
Then $\pi_\wedge(\eta)= 2\omega$, where $\omega$ denotes the
invariant 2-form of Proposition \ref{cycle}.
\end{prop}
\begin{proof}
Again let $n:=1/2$, and set $\tilde P_1:= \fett{t}^{n*}_{n}\,\fett{t}^{n}_{n}$.
From \eqref{ft}, it follows that
$f\rt \tilde P_1=\sigma_{n}(S(f_{(1)})) \tilde P_1\hs \sigma_{n}(f_{(2)})$
for $f\in\su$. By
\eqref{rhoo}, $\rho_{n}^\circ(K^2)= \sigma_{n}(K^{-2})$.
Thus, with $\eta$ given by \eqref{etat},
\begin{align*}
f\rt \pi_\wedge(\eta)
&= \tr\, \sigma_{n}(K^{-2}) \sigma_{n}(S(f_{(1)}))
(1-2\tilde P_1)\dd \tilde P_1\wedge \dd \tilde P_1
\sigma_{n}(f_{(2)})\\
&= \tr\, \sigma_{n}( f_{(2)}  K^{-2}S(f_{(1)}))
(1-2\tilde P_1)\dd \tilde P_1\wedge \dd \tilde P_1
= \vare(f)\, \pi_\wedge(\eta),
\end{align*}
where we used the cyclicity of the trace and $K^2 fK^{-2} = S^{2}(f)$ for
$f\in\su$. Hence $\pi_\wedge(\eta)$ is an invariant 2-form.
By Proposition \ref{cycle},
$\omega$ is also invariant and $\Omega^{\wedge 2}\cong \PS \omega$. Therefore
$\pi_\wedge(\eta)= c\hs \omega$, where $c=\int \pi_\wedge(\eta)$ and, by \eqref{cyc},
\begin{align*}
\int\!\hsp \pi_\wedge(\eta)
=h(q\tr_{\theta^{-1}}(1\!-\!2\tilde P_1)\dif_E(\tilde P_1)\dif_F(\tilde P_1)
\hsp-\hsp q^{-1}\tr_{\theta^{-1}}(1\!-\!2\tilde P_1)\dif_F(\tilde P_1)\dif_E(\tilde P_1)).
\end{align*}
The proof of Proposition \ref{ori} shows that
$\tr_{\theta^{-1}}(1-2\tilde P_1)\dif_F(\tilde P_1)\dif_E(\tilde P_1)=-q$ and
$\tr_{\theta^{-1}}(1-2\tilde P_1)\dif_E(\tilde P_1)\dif_F(\tilde P_1)=q^{-1}$.
Hence $c=2$.
\end{proof}

%*********************************************************
\section{Index computation and Poincar\'e duality}                      \label{IND}
%*********************************************************

The $K$-theoretic version of Poincar\'e duality states that
the additive pairing on $K_*(\CPS)$ determined by the
index map of $D$ is non-degenerate \cite{C1,GFV}.
Since $K_1(\CPS)=0$ \cite{MNW}, it  suffices, in our case, to verify the
non-degeneracy of the pairing
$\ip{\cdot}{\cdot}_D: K_0(\CPS)\times K_0(\CPS)\ra \Z$,
$$
\ip{[P]}{[Q]}_D:=\ind((P\otimes JQJ^*)\,\bar\dif_F\,(P\otimes JQJ^*)),
$$
where $P\in \M_{n\times n}(\CPS)$ and $Q\in \M_{m\times m}(\CPS)$ are
projections, $\bar\dif_F$ denotes the lower-left entry of $D$,
and  $(P\otimes JQJ^*)\,\bar\dif_F\,(P\otimes JQJ^*)$ is an unbounded
Fredholm operator mapping from its domain in
$(P\otimes JQJ^*)\ov{M}_{-1}^{n\times m}$ into the Hilbert space
$(P\otimes JQJ^*)\ov{M}_{1}^{n\times m}$.

Recall that
an unbounded Fredholm operator $F$ is an operator between Hilbert spaces
with dense domain, finite-dimensional kernel, and
finite-co\-di\-men\-sional range.
Its index $\ind(F)$ is the difference between
the dimensions of kernel and cokernel or, equivalently,
$\ind(F)=\dim(\ker F)- \dim(\ker F^*)$.
For an $\su$-equivariant Fredholm operator $F$ (i.e.,
$F$ commutes with the action of $\su$ on the Hilbert space),
kernel and cokernel carry a finite-dimensional representation
of $\su$ and we define
$$
\qind(F):=\tr_{\ker F}K^2 - \tr_{\ker F^*}K^2.
$$

The next lemma is the key to index computations of $D$.

\begin{lem}                                                 \label{ker}
For $N\in \Z$, let  $\bar\dif_E^N$ and  $\bar\dif_F^N$
denote the closures of the
operators $\dif_E : M_N\ra \ov M_{N-2}$ and
$\dif_F : M_N\ra \ov M_{N+2}$, respectively, given by \eqref{dif}.
Then
\begin{align*}
&\ker(\bar\dif_E^{N})=\{0\},\quad
\ker(\bar\dif_F^{N})=\Lin\{ v^{N/2}_{N/2,j}\,:\,
j=-\mbox{$\frac{N}{2}$},\ldots,\mbox{$\frac{N}{2}$}\},& & N>0, \\
&\ker(\bar\dif_E^{N})=\Lin\{ v^{|N|/2}_{N/2,j}\,:\,
j=-\mbox{$\frac{|N|}{2}$},\ldots,\mbox{$\frac{|N|}{2}$}\},\quad
\ker(\bar\dif_F^{N})=\{0\},& &N<0,
\end{align*}
and $\ker(\bar\dif_E^{0})=\ker(\bar\dif_F^{0})=\C v^0_{00}$.
\end{lem}
\begin{proof}
We give the  proof for $\bar\dif_F^{N}$, the other case is similar.
By \eqref{dif},
$v^{l}_{N/2,j}$ is an eigenvector of  $|\bar\dif_F^{N}|$ with corresponding
eigenvalue $\alpha^l_{N/2}$.  The claim follows
from $\ker(\bar\dif_F^{N})=\ker(|\bar\dif_F^{N}|)$ since
$\alpha^l_{N/2}=0$ if and only if  $N\geq 0$ and $l=N/2$.
\end{proof}

We turn now to an explicit computation of indices.

\begin{prop}                                            \label{compi}
For $N\in\Z$,
$$
\ind(P_N\hs\bar\dif_F\hs P_N)=N,\quad
\qind(P_N\hs \bar\dif_F\hs P_N)=[N]_q.
$$
\end{prop}
\begin{rem}
Note that $\ip{[\frac{1}{2}\tau]}{[P_N]}=\qind(P_N\hs\bar\dif_F\hs P_N)=\ip{[\phi]}{[P_N]}$.
The first equality already been established in \cite{NT} from a local index formula.
Whether $\ip{[\phi]}{[P_N]}$ can be computed from a local index formula
will be discussed elsewhere.
\end{rem}
\begin{proof}
Let $n:=N/2>0$ and set $\tilde P_N:=\fett{t}^{n*}_{n}\,\fett{t}^{n}_{n} $.
By the definition of $P_N$ in \eqref{PN},
$\ker(P_N\hs\bar\dif_F\hs P_N)
=\rho_{n}^\circ(K)^{-1}\ker(\tilde P_N\hs\bar\dif_F\hs \tilde P_N)$.
Suppose we are given an  $\fett{x}_{-1}\in \ov M_{-1}^{2n+1}$  such that
$\tilde P_N\hs\bar\dif_F\hs \tilde P_N\fett{x}_{-1}=0$. From
$\dif_F(\fett{t}^{n*}_{n})\sim \fett{t}^{n*}_{n-1}$,
$\fett{t}^{n}_{n}\,\fett{t}^{n*}_{n-1}=0$, $\fett{t}^{n}_{n}\,\fett{t}^{n*}_{n}=1$
and Equation \eqref{difab},
we conclude that $\tilde P_N\hs\bar\dif_F\hs \tilde P_N\fett{x}_{-1}=0$ if and
only if $\fett{t}^{n*}_{n}\hs\bar\dif_F^{N-1}\hs (\fett{t}^{n}_{n}\fett{x}_{-1})=0$.
From \cite[Lemma 6.5(iii)]{SW2}, it follows that $\ker({t^{n*}_{n,-n}})=0$.
Thus
$\fett{t}^{n*}_{n}\hs\bar\dif_F^{N-1}\hs (\fett{t}^{n}_{n}\fett{x}_{-1})\hsp=\hsp0$ if and
only if $\bar\dif_F^{N-1}\hs (\fett{t}^{n}_{n}\fett{x}_{-1})\hsp=\hsp0$
since $\fett{t}^{n*}_{n}\hsp=\hsp({t^{n*}_{n,-n}},\ldots, {t^{n*}_{n,n}})^\trans$.
As $\fett{t}^{n}_{n}\fett{x}_{-1}\in\ov M_{N-1}$, Lemma \ref{ker}  shows that
$$
\ker(\tilde P_N\hs\bar\dif_F\hs \tilde P_N)
= \{\fett{x}_{-1}\in\ov M_{-1}^{2n+1}\,:\, \fett{t}^{n}_{n}\fett{x}_{-1}
\in \Lin\{ v^{n-1}_{n-1,-n+1},\ldots,v^{n-1}_{n-1,n-1}\}\hs\}.
$$
By Lemma  \ref{MN}, we can find $\fett{x}_{-1}^{j}\in M_{-1}^{2n+1}$
%$\fett{x}_{-1}^{-(n-1)},\ldots, \fett{x}_{-1}^{n-1}\in M_{-1}^{2n+1}$
such that
$\fett{t}^{n}_{n}\fett{x}_{-1}^j=v^{n-1}_{n-1,j}$, $j=-(n-1),\ldots,n-1$. Set
$\tilde{\fett{x}}_{-1}^j:=\rho_{n}^\circ(K)^{-1}(\fett{x}_{-1}^j)$.
By the preceding, $\tilde{\fett{x}}_{-1}^j\in\ker(P_N\hs\bar\dif_F\hs P_N)$.
Suppose that
$\tilde{\fett{y}}_{-1}^{-n},\ldots, \tilde{\fett{y}}_{-1}^{n}$ is another such set.
Then
$$
P_N(\tilde{\fett{x}}_{-1}^j-\tilde{\fett{y}}_{-1}^j)\hsp=\hsp\rho_{n}^\circ(K^{-1})
\fett{t}^{n*}_{n}\hs\fett{t}^{n}_{n}(\fett{x}_{-1}^j-\fett{y}_{-1}^j)
\hsp=\hsp\rho_{n}^\circ(K^{-1})\fett{t}^{n*}_{n}\hs(v^{n-1}_{n-1,j}-v^{n-1}_{n-1,j})\hsp=\hsp0,
$$
so $P_N\hs\tilde{\fett{x}}_{-1}^j=P_N\hs\tilde{\fett{y}}_{-1}^j$ in $P_N\ov M_{-1}^{2n+1}$.
Hence
\begin{equation}                                             \label{kerF}
\ker(P_N\hs\bar\dif_F\hs P_N)\cong \Lin\{ v^{n-1}_{n-1,-(n-1)},\ldots,v^{n-1}_{n-1,n-1}\}.
\end{equation}
Now we consider $(P_N\hs\bar\dif_F\hs P_N)^*=P_N\hs\bar\dif_E\hs P_N$.
Using $\dif_E(\fett{t}^{n*}_{n})=0$, the same reasoning as above shows that
$\tilde P_N\hs\bar\dif_E\hs \tilde P_N\fett{x}_{1}=0$ for
$\fett{x}_{1}\in \ov M_{1}^{2n+1}$ if and only if
 $\bar\dif_E\hs (\fett{t}^{n}_{n}\fett{x}_{1})
=\bar\dif_E^{N+1}\hs (\fett{t}^{n}_{n}\fett{x}_{1})=0$. Lemma \ref{ker} implies
$\fett{t}^{n}_{n}\fett{x}_{1}=0$, thus
$P_N (\rho_{n}^\circ(K)^{-1}\fett{x}_{1})=0$ and therefore
$\ker(P_N\hs\bar\dif_E\hs P_N)=\{0\}$. Summarizing, we get
\begin{align*}
\ind(P_N\hs\bar\dif_F\hs P_N)=\dim(\ker(P_N\hs\bar\dif_F\hs P_N))=2(n-1)+1=N,\\
\qind(P_N\hs\bar\dif_F\hs P_N)=\tr_{\ker(P_N\hs\bar\dif_F\hs P_N)}K^2
=\mbox{$\sum$}_{j=-n+1}^{n-1}\hs q^{2j}=[N]_q
\end{align*}
by \eqref{kerF}. This proves the proposition for $N>0$.
The case $N=0$ is trivial, and the case $N<0$ is proved similarly with the
role of $\bar\dif_F$ and $\bar\dif_E$ interchanged.
\end{proof}

It has been shown in \cite{MNW} that $\CPS\cong \C\hs 1 + \cK(\lN)$, where
$\cK(\lN)$ denotes the compact operators on $\lN$, and
$\KO(\CPS)\cong \Z \oplus \Z$. The generator $[(1,0)]$ was taken to be
the identity $P_0=1$, and the other generator  $[(0,1)]$ the 1-dimensional
projection onto the first basis vector. From  \cite{MNW} and \cite{Haj},
we conclude that $[P_1]=[(1,1)]$. In particular, $[P_0]$ and $[P_1]$
also generate $\KO(\CPS)$.

Using the results from the previous proposition,
we can easily establish Poincar\'e duality:
\begin{prop}
The pairing $\ip{\cdot}{\cdot}_D: K_0(\CPS)\times K_0(\CPS)\ra \Z$ is given by
$$
\ip{[(k,l)]}{[(m,n)]}_D = km-ln, \quad k,l,m,n\in\Z.
$$
In particular, it is non-degenerate, so Poincar\'e duality holds.
\end{prop}
\begin{proof}
We first show that $\ip{\cdot}{\cdot}_D$ is antisymmetric.
Let $P$ and $Q$ be projection matrices with entries in $\CPS$.
Note that
$\ind((P\otimes Q)\,\bar\dif_F\,(P\otimes Q))=
\ind((Q\otimes P)\,\bar\dif_F\,(Q\otimes P))$ since the flip of
tensor factors is an unitary operation which commutes with the
component wise action of $\bar\dif_F$. Next,
$J^*PJ=JPJ^*$ since $J^2=-1$.
Recall from Section \ref{DO} that $D$ and $J$ are odd operators, i.e.,
$\gamma D=-D\gamma$ and  $\gamma J=-J\gamma$, and $JD=DJ$.
Hence $J^* \bar\dif_F J= \bar\dif_E=\bar\dif_F^*$.
Moreover, $\ind(F)=-\ind(F^*)$ for any Fredholm operator $F$. Thus
\begin{align*}
\ip{[P]}{[Q]}_D &=
\ind((P\otimes JQJ^*)\,\bar\dif_F\,(P\otimes JQJ^*))\\
&=\ind(J( J^*PJ\otimes Q)\,J^*\bar\dif_FJ\,(J^*PJ\otimes Q)J^*)\\
&=-\ind((Q\otimes JPJ^*)\,\bar\dif_F\,(Q\otimes JPJ^*))= - \ip{[Q]}{[P]}_D.
\end{align*}

By the additivity of the pairing, it suffices to consider the generators
$[P_1]=[(1,1)]$ and $[P_0]=[(1,0)]$. From Proposition \ref{compi}, we get immediately
$\ip{[(1,1)]}{[(1,0)]}_D=1$, and the antisymmetry implies $\ip{[(1,0)]}{[(1,1)]}_D=-1$,
$\ip{[(1,0)]}{[(1,0)]}_D=\ip{[(1,1)]}{[(1,1)]}_D=0$.
For $k,l,m,n\in\Z$, we obtain now
\begin{align*}
\ip{[(k,l)]}{[(m,n)]}_D
&=\ip{\,k[(1,1)]+(l\!-\!k)[(1,0)]\,}{\,m[(1,1)]+(n\!-\!m)[(1,0)]\,}_D\\
&=kn-lm,
\end{align*}
which completes the proof.
\end{proof}
%
%

%**************************
\section*{Acknowledgments}
%**************************
%
The author thanks Ludwik D\c{a}browski, Francesco D'Andrea, Giovanni Landi,
and Adam Rennie for stimulating discussions on the subject, and is
especially grateful to an anonymous referee for his very useful comments.
This work was
supported by the DFG fellowship WA 1698/2-1 and the European Commission grant
MTK-CT-2004-509794.

\end{document}